\documentclass[journal,
    ]{IEEEtran}%
\pdfoutput=1

\usepackage{microtype}
\usepackage{graphicx}
\usepackage{subfigure}
\usepackage{booktabs} 
\usepackage[binary-units=true]{siunitx}
\usepackage[
	colorlinks=true,
	hypertexnames=false
	]{hyperref}
\PassOptionsToPackage{linktocpage}{hyperref} 

\hypersetup{
  colorlinks   = true, 
  urlcolor     = magenta, 
  linkcolor    = blue, 
  citecolor   = teal 
}

\usepackage{amsmath,amsfonts,amssymb,amsthm,bbm,graphicx,enumerate,times, mathdots,braket}
\usepackage{mathtools}
\usepackage{tikz}
\usepackage{graphicx}
\usepackage[normalem]{ulem}
\usepackage{float}
\usepackage{comment}
\usepackage{csquotes}
\usepackage{url}
\usepackage{enumerate}  
\usepackage{natbib} 
\renewcommand\cite{\citep}
\usepackage{algorithm}
\usepackage{algorithmic}
\usetikzlibrary{decorations.pathreplacing}

\newtheorem{theorem}{Theorem}

\newtheorem{definition}[theorem]{Definition}

\newtheorem{problem}[theorem]{Problem}
\newtheorem{rem}[theorem]{Observation}

\definecolor{jens}{rgb}{.2,.7,.9}
\definecolor{ingo}{rgb}{.9,.5,.2}

\newcommand{\tu}{Institute of Mathematics, Technische Universit{\"a}t Berlin, Germany}
\newcommand{\fu}{Dahlem Center for Complex Quantum Systems, Freie Universit{\"a}t Berlin, Germany}

\title{Tensor network approaches for learning non-linear dynamical laws}

\author{
    \IEEEauthorblockN{
        A.~Goe{\ss}mann\IEEEauthorrefmark{1}\IEEEauthorrefmark{4},
        M.~G{\"o}tte\IEEEauthorrefmark{1}\IEEEauthorrefmark{4},
        I.~Roth\IEEEauthorrefmark{3},
        R.~Sweke\IEEEauthorrefmark{3},
        G.~Kutyniok\IEEEauthorrefmark{4},
        and J.~Eisert\IEEEauthorrefmark{3}}

    \IEEEauthorblockA{\IEEEauthorrefmark{4}\tu{}}

    \IEEEauthorblockA{\IEEEauthorrefmark{3}\fu{}}

    \thanks{\IEEEauthorrefmark{1} Both authors contributed equally. (goessmann@tu-berlin.de)
    }
}

\makeatletter 
  \hypersetup{pdftitle = {Tensor network approaches for learning non-lineardynamical laws},
         pdfauthor = {
            Alex Goessmann, 
            Michael Goette, 
            Ingo Roth, 
            Ryan Sweke, 
            Gita Kutyniok, 
            Jens Eisert
         },
         pdfsubject = {
         non-linear dynamics, 
         automate theory building,
         tensor networks,
         machine learning
         },
         pdfkeywords = {
         tensor network,
         machine learning,
         dynamical laws,
         theory building, 
         partial differential equation, 
         tensor train,
         matrix product state, 
         rank-adaptivity, 
         SALSA, 
         alternating least square, 
         governing equation, 
         sparse identification, 
         compressed sensing, 
         curse of dimensionality, 
         multivariate, 
         univariate,
         locality
            }
        }
\makeatother

\begin{document}

\maketitle

\begin{abstract}

Given observations of a physical system, identifying the underlying non-linear governing equation is a fundamental task, necessary both for gaining understanding and generating deterministic future predictions. Of most practical relevance are automated approaches to theory building that scale efficiently for complex systems with many degrees of freedom. To date, available scalable methods aim at a data-driven interpolation, without exploiting or offering insight into fundamental underlying physical principles, such as locality of interactions. In this work, we show that various physical constraints can be captured via tensor network based parameterizations for the governing equation, which naturally ensures scalability. In addition to providing analytic results motivating the use of such models for realistic physical systems, we demonstrate that efficient rank-adaptive optimization algorithms can be used to learn optimal tensor network models without requiring a~priori knowledge of the exact tensor ranks. As such, we provide a physics-informed approach to recovering structured dynamical laws from data, which adaptively balances the need for expressivity and scalability.
\end{abstract}

\section{Introduction}

A core method in the natural sciences is the inference of an equation governing the dynamics of a physical system from observations. Finding a \emph{governing equation} typically provides an explanatory understanding of the dynamics that goes significantly beyond the mere identification of statistical relations and patterns in the observed data.
While modern computing technologies have facilitated automated access to vast amounts of data, from a wide variety of dynamical systems, flexible and scalable approaches to automatically extract governing equations from this data are rare.
For a traditional theorist, the process of inferring a suitable governing equation is
guided by expert intuition and domain knowledge, which allows one to exploit physically motivated constraints, such as admissible correlation structures.
Automated approaches to this task have to balance the need for expressivity, which allows for the exploration of large theory spaces, with the need for scalability, which facilitates application to systems with many degrees of freedom.
This is exactly the realm where the `traditional approach' is bound to fail and automated approaches may turn out most valuable for consistent theory building. 

To be more precise, let us model the state of a dynamical system by $d$ real variables $(x_1,\dots,x_d)=:x$ and its time evolution by smooth trajectories $t \mapsto x(t)\in \mathbb{R}^d$.
The dynamical laws of many physical systems are formulated by a suitable differential operator $\mathcal D$, such as the derivatives $\frac{d}{dt}$ or $\frac{d^2}{dt^2}$, acting on the trajectory coordinatewise as
\begin{equation}
\nonumber
    \mathcal{D} x(t)=
    \left[ \begin{array}{cc}
          \mathcal{D} x_1(t)   \\
         \mathcal{D} x_2(t)   \\
         \, \, \vdots  \\
         \mathcal{D} x_d(t)
    \end{array}  \right]      
        =: \left[ \begin{array}{c}
         f_1(x(t))   \\
         f_2(x(t))   \\
         \quad \vdots  \\
         f_d(x(t))
    \end{array}  \right] =
    f\big( x(t) \big)\, .
\end{equation}
To learn the governing equation of a system, we assume access to potentially noisy estimates of different states $x^j$ and the corresponding evaluation of the differential operator or, equivalently, $f(x^j)$. 
Such data can, for example, be obtained from time-series data by finite difference approximations of $\mathcal{D}$. Given a \textit{hypothesis set} $\mathcal{F}$ of functions $\tilde{f}: \mathbb{R}^d \rightarrow \mathbb{R}$ on the state space, which represents the expert's intuition and the domain knowledge as a prior assumption, we formulate the recovery problem in the following fashion:

\begin{problem}[Governing equation recovery] \label{prob:GoverningEquation}
    Identify the governing equation $f=[f_1\dots f_d]$ with $f_k\in \mathcal{F}$ from the given observations $\{x^j,y^j:=f(x^j)\}_{j=1}^m$.
\end{problem}

The identification of a function from evaluations at sample points is the central task of \emph{supervised machine learning} \cite{bishop_pattern_2006}. Note, however, that traditional machine learning methods typically aim to interpolate complex relations. Theory building, in contrast, is
more ambitious and is interested in an exact recovery of an 
\emph{interpretable governing equation}, 
which motivates our approach.

Simple hypothesis sets $\mathcal F$ are linear spaces spanned by 
\emph{basis functions} $\{\phi_1(x),\dots,\phi_p(x)\}$. In this space, functions $f_l$ are represented by their linear coefficients $\theta_{il}$ and function evaluations are modeled by multiplication with a dictionary matrix $\Phi_{ji}:=\phi_i(x^j)$ as
\begin{equation}
    y_l^j=\sum_{i=1}^p \phi_i(x^j) \theta_{il} \,.
\label{DefCoefficientRecovery}
\end{equation}
We employ a graphical notation%
\footnote{As is commonly done in this context, we sketch tensors of different orders by rectangles with legs representing its indices.
Connections of legs to different tensor indicates scalar products performed in the respective index spaces, see App.~\ref{app:notation}.},
in which \eqref{DefCoefficientRecovery} reads
\begin{center}
\begin{tikzpicture}[scale=0.3]
    \draw (-6,0)--(-4.5,0) node[midway,above] {$j$};
    \draw (-4.5,-1) rectangle (-2.5,1);
    \coordinate[label=above:$y$] (A) at (-3.5,-0.8);
    \draw (-2.5,0)--(-1,0) node[midway,above] {$l$};
    \coordinate[label=above:{$=$}] (A) at (0,-0.8);
    \draw (1,0)--(2.5,0) node[midway,above] {$j$};
    \draw (2.5,-1) rectangle (4.5,1);
    \coordinate[label=above:{$\Phi$}] (A) at (3.5,-0.8);
    \draw (4.5,0)--(6,0) node[midway,above] {$i$};
    \draw (6,-1) rectangle (8,1);
    \coordinate[label=above:$\theta$] (A) at (7,-0.8);
    \draw (8,0)--(9.5,0) node[midway,above] {$l$};
\end{tikzpicture}.
\end{center} 

Problem~\ref{prob:GoverningEquation} then amounts to the linear inverse problem of recovering the coefficients $\theta_{il}$ representing the system $[f_1\dots f_d]$.
Since many dynamical systems need only a few dominant elementary functions to represent their governing equation, Ref.~\cite{brunton_discovering_2016} proposes the \emph{sparse identification of non-linear dynamics (SINDy)}.
A sparsity assumption on the coefficient vector effectively amounts to restricting $\mathcal F$ to the union of its lower-dimensional sub-spaces.
In addition, it can resolve the ill-positioning of the linear inverse problem in the regime where $\Phi$ has a non-trivial kernel, e.g.\ due to insufficient data \cite{schaeffer_extracting_2018}. 
The assumption might also be interpreted as an implementation of the principle of \textit{Occam's Razor}.

For systems with many degrees of freedom suitable linear function dictionaries for multivariate functions become inadmissibly large.
This `curse of dimensionality' yields a severe limitation for the scalability of SINDy to the regime where automated theory building would be most helpful.
Similarly, approaches using \textit{symbolic regression} \cite{schmidt_distilling_2009, kusner_grammar_2017-1,ouyang_sisso:_2018, li_neural-guided_2019} for automated theory building will generate exponentially many combinations of basis functions, again severely limiting their applicability for systems with many degrees of freedom. Furthermore, the optimization of non-linear coefficients in these functions sets leads to infinite dimensional linear hulls. This problem can be mitigated by \emph{multi-linear} parameterization schemes \cite{cohen_expressive_2015,stoudenmire_supervised_2016,levine_deep_2017,gels_multidimensional_2019}. 
At the heart of these parameterization schemes is the insight that multivariate function dictionaries often feature a tensor structure. 

In this work, building on this insight, we show that fundamental physical principles `naturally' allow us to break the curse of dimensionality. 
By taking into consideration physical principles such as locality and symmetry constraints, we derive tensor network formats which allow for the exploitation of expert knowledge about a system in the data-driven governing equation recovery process.
Specifically, we first argue that certain tensor networks, such as low-rank tensor trains, ensure a constrained correlation structure of the coefficients for each function component $f_l$. 
Secondly, we show how to incorporate these tensor networks into larger networks that ensure additional assumptions about the correlations between the component functions of the governing equations.
We substantiate our findings by studying classes of one-dimensional systems with local interactions, where we explicitly derive the ranks required for parameterization in the respective formats. As example systems in such classes we discuss variants of the Fermi-Pasta-Ulam-Tsingou model.
We introduce state-of-the-art rank-adaptive numerical optimization schemes for the recovery of governing equations of unknown rank, parameterized by the previously introduced tensor network models.
Importantly, the rank-adaptivity of our approach ensures that our assumptions on the correlation structure are not a hard-coded restriction on the explored function space, thus we do not require detailed knowledge about suitable ranks and tensor structures for a specific problem.
Our approach rather uses structural insights to identify numerically feasible corners from which the problem can be explored, adaptively increasing the computational effort. 
This work therefore provides a general framework in which traditional concepts of theory building are combined with data-driven optimization algorithms, facilitating progress towards efficient and broadly applicable automated theory discovery.

The remainder of this work is structured as follows: 
After introduction of a tensor parameterization for multivariate functions in Section~\ref{sec:multivariate}, we develop tensor network models suitable for the coefficient vector in Section~\ref{sec:models}. 
Section~\ref{sec:optimization} then introduces algorithms based on alternating least squares optimization and its rank-adaptive generalization, which will be tested in specific recovery tasks in Section~\ref{sec:examples}.

\section{Parameterization of multivariate functions} \label{sec:multivariate}

\begin{figure}
    \centering
\begin{tikzpicture}[scale=0.3,thick]

\begin{scope}[shift={(0,-0.75)}]
    \draw (-15,-5.5) rectangle (-13,-3.5);
    \draw (-15,-4.5)--(-17,-4.5) node[midway,above] {$j$};  
    \coordinate[label=above:$y$] (A) at (-14,-5.5);
    \coordinate[label=above: {\Large =}] (A) at (-11,-5.35);
\end{scope}

\begin{scope}[yscale=-1,xscale=-1,shift={(0,9)}]

    \draw (-1,-1) rectangle (3.5,1);
     \coordinate[label=below:$\psi_{i_1}(x^{\, j}_1)$] (A) at (1.25,-1.3);
     \draw (-1,0)--(-2.5,0) node[midway,above] {$i_1$};

    \begin{scope}[shift={(0,-3)}]  
       \draw (-1,-1) rectangle (3.5,1);
        \coordinate[label=below:$\psi_{i_2}(x^{\, j}_2)$] (A) at (1.25,-1.3);
        \draw (-1,0)--(-2.5,0) node[midway,above] {$i_2$};
    \end{scope}

    \coordinate[label=below:$\vdots$] (A) at (1.25,-7);
    \coordinate[label=below:$\vdots$] (A) at (5,-7);

    \begin{scope}[shift={(0,-7.5)}]  
        \draw (-1,-1) rectangle (3.5,1);
        \coordinate[label=below:$\psi_{i_d}(x^{\, j}_d)$] (A) at (1.25,-1.3);
        \draw (-1,0)--(-2.5,0) node[midway,above] {$i_d$};
    \end{scope}
 
    \draw[fill=black] (7,-3.75) circle (0.25cm);
     \coordinate[label=right:$ \delta$] (A) at (7,-4.25);
    \draw (7,-3.75) to[bend right=50] (3.5,0);
    \draw (7,-3.75) to[bend right=10] (3.5,-3);
    \draw (7,-3.75) to[bend left=50] (3.5,-7.5);
    \draw (7,-3.75)--(9,-3.75) node[midway,above] {$j$};   

    \draw (-4.5,1) rectangle (-2.5,-8.5); 
    \coordinate[label=below:$\theta$] (A) at (-3.5,-4.75);
   
    \end{scope}
    
     \draw[dashed] (2.25,-10.5)--(-8.5,-10.5) node[midway, below, yshift=0pt,]{Dictionary matrix $\Phi$}--(-8.5,0)--(2.25,0)--(2.25,-10.5);
\end{tikzpicture}
    \caption{Tensor recovery as a linear inverse problem with observations $y$ generated by contraction of the coefficient tensor $\theta$ on the dictionary matrix $\Phi$ (dashed). Employing a Hadamard product, represented by the contraction with the tensor $\delta$, $\Phi$ is decomposed into univariate dictionary matrices $\psi_{i_k}(x_k^j)$.}
    \label{fig:tensor_recovery}
\end{figure}
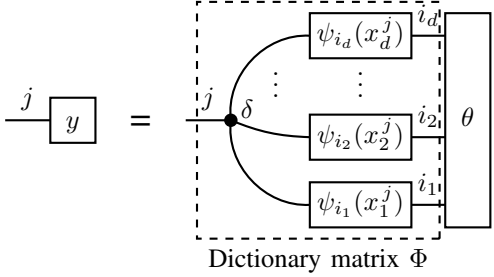

Let us first discuss the natural method to build multivariate function spaces as tensor products of univariate function spaces, which will then enable us to find suitable hypothesis subsets $\mathcal{F}$ for Problem \ref{prob:GoverningEquation}.
To this end, 
we assume a set of linearly independent basis functions $\{\psi_{i}:\mathbb{R}\to\mathbb{R}\}_{i=1}^{\tilde p}$ that act on the individual coordinates and take the products 
\begin{equation}
    \phi_{i_1\ldots i_d}(x_1,\dots,x_d):=\psi_{i_1}(x_1)\cdot\psi_{i_2}(x_2)\cdot \ldots \cdot \psi_{i_d}(x_d)\,.
\label{DefTensorBasis}
\end{equation}
Note that our approach in principle also allows for choosing different function sets for each coordinate. 
In graphical notation, the coordinates $x_k$ are represented by vectors $\psi(x_k)$ of dimension $\tilde{p}$, which build the tensor representation $\phi(x_1,\dots,x_d)$ of the whole state by an outer product as
\begin{center}
\begin{tikzpicture}[scale=0.3]
\begin{scope}[shift={(1,0)}]
    \draw (-8,-1) rectangle (-1,1);
    \coordinate[label=above:{$\phi(x_1,...,x_d)$}] (A) at (-4.5,-0.9);
    \draw (-7.25,-1)--(-7.25,-2.25) node[midway,left] {$i_1$};
    \draw (-5.75,-1)--(-5.75,-2.25) node[midway,left] {$i_2$};
    \draw (-1.75,-1)--(-1.75,-2.25) node[midway,left] {$i_d$};
    \coordinate[label=below:{$\dots$}] (A) at (-4.25,-1.2);
\end{scope}    
    \coordinate[label=above:{$=$}] (A) at (1.25,-0.6);
    
    \draw (2.5,-1) rectangle (5.5,1);
    \coordinate[label=above:{$\psi(x_1)$}] (A) at (4,-0.9);
    \draw (4,-1)--(4,-2.25) node[midway,left] {$i_1$};
    
    \draw (6.5,-1) rectangle (9.5,1);
    \coordinate[label=above:{$\psi(x_2)$}] (A) at (8,-0.9);
    \draw (8,-1)--(8,-2.25) node[midway,left] {$i_2$};
    
    \coordinate[label=above:{$\dots$}] (A) at (11,-0.4);
    
    \draw (12.5,-1) rectangle (15.5,1);
    \coordinate[label=above:{$\psi(x_d)$}] (A) at (14,-0.9);
    \draw (14,-1)--(14,-2.25) node[midway,left] {$i_d$};
\end{tikzpicture}
\end{center} 
This tensor representation is referred to as the \emph{coordinate major} by \citet{gels_multidimensional_2019}. The number $p=\tilde{p}^d$ of basis functions, constructed by this product ansatz, grows exponentially in the number $d$ of variables. This gives such function spaces a high expressivity but also manifests the \emph{curse of dimensionality} that linear methods, such as SINDy, suffer for high-dimensional systems.
At the same time, the tensor product structure (\ref{DefTensorBasis}) allows for calculation and storage of the dictionary matrix $\Phi$ with linear demand in the number $d$ of variables, by evaluating the univariate functions $\psi_{i_k}$ at the respective coordinates $x_k^j$ of the observed states. Storing the function evaluations in the univariate dictionary matrices $\Psi^k:=\psi_{i_k}(x_k^j)$ unravels the structure of the multivariate dictionary matrix $\Phi$ as a Hadamard product of $\Psi^k$ along the data index $j$, which is represented in Fig.~\ref{fig:tensor_recovery} by contraction with a delta tensor $\delta$ with unit entries only on its hyperdiagonal. 

Each function in the span of the products \eqref{DefTensorBasis} is represented by its \textit{coefficient tensor} $\theta$ through its basis decomposition
\begin{equation}
    f^{\theta}(x_1,\dots,x_d):= \sum_{i_1,\dots,i_d=1}^{\tilde{p}} \theta_{i_1\dots i_d} \psi_{i_1}(x_1)\cdot\ldots\cdot \psi_{i_d}(x_d) \,.
\label{eq:TensorParameterization}    
\end{equation}

If the univariate basis functions $\{\psi_{i}\}$ are orthonormal, the euclidean scalar product in the coefficient tensor space equals the scalar product of the indexed functions, which is induced by the scalar products of the univariate functions. In contrast to the sparsity hypothesis exploited in the SINDy method \cite{ho_recovery_2018}, we will in this work derive hypothesis sets that are invariant under span-preserving transformations of the basis functions. We will thus always be able to parameterize with respect to orthonormal basis functions, which ensures the equivalence of the scalar products on the coefficient space and the function spaces.

The evaluation of the functions \eqref{eq:TensorParameterization} at the states $x^j$ amounts to a contraction of the dictionary matrix $\Phi$ with the corresponding coefficient tensors. 
Identifying for each function $f_l$ in Problem~\ref{prob:GoverningEquation} the corresponding coefficient tensor $\theta_l$ becomes therefore a multi-linear inverse problem. 
Our strategy to avoid the `curse of dimensionality' is to restrict the exponentially large tensor space to a tractable hypothesis subset $\Theta \subset \mathbb{R}^{p}$, for which we pose the following problem:

\begin{problem}[Tensor recovery]
    Given a dictionary matrix $\Phi$ and observations  $y=\Phi \theta$,
    recover the tensor $\theta$ under the assumption $\theta \in \Theta$ with hypothesis set $\Theta \subset \mathbb{R}^{p}$.
    \label{prob:TensorRecovery}
\end{problem}

One approach to this problem proposed by \citet{klus_tensor-based_2019} takes the perspective of kernel methods: The function \eqref{DefTensorBasis} can be regarded as a feature map lifting the states $x$ into a tensor space $\mathbb{R}^{p}$ with coordinate maps $\phi_{i_1,\ldots,i_d}$.
Then, inspired by the \emph{representer theorem} 
\cite{hofmann_kernel_2008} 
the subspace spanned by the represented states $x^j$ can be taken as the hypothesis set $\Theta$.
The \emph{MANDy method} originally introduced in Ref.~\cite{gels_multidimensional_2019} efficiently identifies the minimal norm solution of the associated least squares problem 
by directly calculating the pseudo-inverse which maps to this hypothesis set.
This approach has the drawback that the hypothesis is chosen by the data and does not directly take into account physical considerations.
Our tensor network approach, which we will introduce in the next section, naturally overcomes this shortcoming by using locality and symmetry as guiding principles.

\section{Tensor network hypothesis manifolds}\label{sec:models}

Problem~\ref{prob:TensorRecovery} naturally leads to the main question of this work:
 \emph{What are natural hypothesis sets $\Theta$?}
In this section we develop tensor network models for the hypothesis set $\Theta$ that are informed by physical paradigms but remain highly adaptive. 

\subsection{Correlations and separation ranks} \label{subsec:separation}
In the context of tensor network representation of many-body quantum states, constraints on the structure of correlations relate directly to the locality of interactions in the physical system and provide a solid theoretical motivation for such representations \cite{Schuch_MPS,eisert_area_2010}.
Motivated by these insights, we define here a notion of correlation, founded on the decomposition properties of functions $f:\mathbb{R}^d\rightarrow\mathbb{R}$ into univariate functions, which represent the variables of a dynamical system. This equips us with a theoretical basis to discuss locality principles for governing equations.

The most extreme case is the absence of correlations:
A multivariate function is called \textit{separable}, if it is the product of univariate functions \cite{beylkin_algorithms_2005}. Functions $f^{\theta}$ parameterized by tensors $\theta$ (see Eq. \eqref{eq:TensorParameterization}) are thus separable, if and only if the parameterizing tensor is elementary, i.e. it is an outer product $v^1\otimes\dots\otimes v^d$ of vectors $v^k\in \mathbb{R}^{\tilde{p}}$. 
While generic tensors do not satisfy such a decomposition, one can always decompose a tensor into a sum of elementary tensors, which is called a $\mathrm{CP}$-decomposition. 
The $\mathrm{CP}$-rank $r(\theta)$ of a tensor $\theta$ is the smallest number of elementary tensors appearing in a $\mathrm{CP}$-decomposition.
Since we have chosen linearly independent basis functions $\psi_{i}$, $r(\theta)$ determines the minimal number of separable functions for a decomposition 
\begin{align}
f^{\theta}(x_1,\dots,x_d)=\sum_{l=1}^{r(\theta)} f^{v^{1,l}}(x_1)\cdot
\dots \cdot f^{v^{d,l}}(x_d), 
\end{align}
where $v^{k,l}\in \mathbb{R}^{\tilde{p}}$ for $1\leq k \leq d$ and $1\leq l \leq r$. 
The $\mathrm{CP}$-rank of a multivariate function provides a measure of correlations with respect to the partition $\{\{x_1\},\ldots,\{x_d\}\}$ of its variables. 
Similarly, given an arbitrary partition $\mathcal{P}$ of the variables $\{x_1,\dots,x_d\}$ into disjoint subsets, we define a separable function with respect to $\mathcal{P}$ as a function which can be written as the product of functions which depend only on variables in a subset of the partition. Accordingly, the separation rank $r_{\mathcal{P}}$ is then defined as the minimal number of separable functions, with respect to $\mathcal{P}$, such that $f^{\theta}$ can be decomposed into a sum of these functions. 
Guided by the conceptual framework of many-body physics \cite{eisert_area_2010}, we understand such separation ranks as a quantification of the correlation structure of a multivariate function.
Note, however, that using $r_{\mathcal{P}}$ directly as a measure for correlation comes with stability issues related to the notion of border rank \cite{BiniLottiRomani:1980}.
While a generic tensor has full separation rank with respect to any partition, many functions arising in the sciences have collections of partitions with low separation ranks.
This can be understood as arising from a notion of locality inherent in the underlying governing equations.
In these cases, it is often possible to obtain a low-rank tensor network decomposition of the tensor $\theta$ encoding the function.

\begin{figure}
    \centering
   \begin{tikzpicture}[scale=0.3,thick]

\begin{scope}[shift={(-15,0)},xscale=-1,yscale=1]
    \coordinate[label=above:$a)$] (A) at (0,0.5);
    \draw (-7.5,1) rectangle (-5.5,-10); 
    \coordinate[label=above:$\theta$] (A) at (-6.5,-5.5);
    \draw (-4,0)--(-5.5,0) node[midway,above] {$i_d$};
    \draw (-4,-1) rectangle (-1,1);
    \coordinate[label=above:$\psi(x_d)$] (A) at (-2.5,-1);
    
    \draw (-4,-5.5)--(-5.5,-5.5) node[midway,above] {$i_2$};
        \draw (-4,-6.5) rectangle (-1,-4.5);
    \coordinate[label=above:$\psi(x_2)$] (A) at (-2.5,-6.5);
    
    \draw (-4,-9)--(-5.5,-9) node[midway,above] {$i_1$};
        \draw (-4,-10) rectangle (-1,-8);
    \coordinate[label=above:$\psi(x_1)$] (A) at (-2.5,-10);
    
    \coordinate[label=above:$\vdots$] (A) at (-4.75,-3.5);
    \coordinate[label=above:$\vdots$] (A) at (-2.5,-3.5);
    
\end{scope}

\begin{scope}[shift={(-6,0)},xscale=-1,yscale=1]
    \coordinate[label=above:$b)$] (A) at (0,0.5);
    \draw (-7.5,-8) rectangle (-5.5,-10); 
    \coordinate[label=above:$\alpha$] (A) at (-6.5,-9.8);
    \draw (-7.5,-6.5) rectangle (-5.5,1); 
    \coordinate[label=above:$\beta$] (A) at (-6.5,-3.55);
     \draw (-6.5,-8)--(-6.5,-6.5) node[midway,right] {$l_1$};  

    \draw (-4,0)--(-5.5,0) node[midway,above] {$i_d$};
    \draw (-4,-1) rectangle (-1,1);
    \coordinate[label=above:$\psi(x_d)$] (A) at (-2.5,-1);
    
    \draw (-4,-5.5)--(-5.5,-5.5) node[midway,above] {$i_2$};
        \draw (-4,-6.5) rectangle (-1,-4.5);
    \coordinate[label=above:$\psi(x_2)$] (A) at (-2.5,-6.5);
    
    \draw (-4,-9)--(-5.5,-9) node[midway,above] {$i_1$};
        \draw (-4,-10) rectangle (-1,-8);
    \coordinate[label=above:$\psi(x_1)$] (A) at (-2.5,-10);
    
    \coordinate[label=above:$\vdots$] (A) at (-4.75,-3.5);
    \coordinate[label=above:$\vdots$] (A) at (-2.5,-3.5);
    
\end{scope}


\begin{scope}[xscale=1,yscale=1,shift={(9.5,0)}]
    \coordinate[label=above:$c)$] (A) at (-6.5,0.5);
\draw (-1,-1) rectangle (1,1);
 \coordinate[label=above:$A^d$] (A) at (0,-1);
 \draw (-1,0)--(-2.5,0) node[midway,above] {$i_d$};
     \draw (-2.5,-1) rectangle (-5.5,1);
    \coordinate[label=above:$\psi(x_d)$] (A) at (-4,-1);
 
  \draw (0,-1)--(0,-2) node[midway,right] {$l_{d-1}$};
   \draw (0,-3.5)--(0,-4.5) node[midway,right] {$l_2$};
    \draw (0,-6.5)--(0,-8) node[midway,right] {$l_1$};
  
\begin{scope}[shift={(0,-5.5)}]  
\draw (-1,-1) rectangle (1,1);
 \coordinate[label=above:$A^2$] (A) at (0,-1);
 \draw (-1,0)--(-2.5,0) node[midway,above] {$i_2$};
      \draw (-2.5,-1) rectangle (-5.5,1);
    \coordinate[label=above:$\psi(x_2)$] (A) at (-4,-1);
\end{scope}

  \coordinate[label=above:$\vdots$] (A) at (0,-3.7);
  \coordinate[label=above:$\vdots$] (A) at (-4,-3.7);  

\begin{scope}[shift={(0,-9)}]  
\draw (-1,-1) rectangle (1,1);
 \coordinate[label=above:$A^1$] (A) at (0,-1);
 \draw (-1,0)--(-2.5,0) node[midway,above] {$i_1$};
      \draw (-2.5,-1) rectangle (-5.5,1);
    \coordinate[label=above:$\psi(x_1)$] (A) at (-4,-1);
\end{scope}

 \end{scope}
\end{tikzpicture}
    \caption{Parameterization of a function by (a) a generic tensor $\theta$, which is in (b) decomposed into tensors $\alpha$ and $\beta$ and in (c) decomposed into the tensor 
    train format by iterative repetition of the decomposition.}
    \label{fig:tensortrain}
\end{figure}
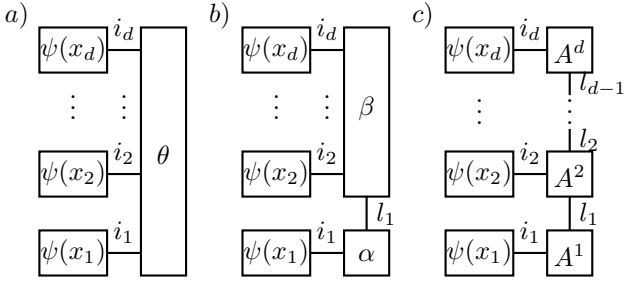

\subsection{Tensor train format}\label{subsec:tt}
A concrete example of an efficient tensor network parameterization with a clear understanding of the underlying locality structure is given by the \emph{tensor train (TT)} format (see Fig.~\ref{fig:tensortrain}c),
referred to as \emph{matrix product state} in the many-body physics literature \cite{quant-ph/0608197}.
Given a one-dimensional collection of coordinates $\{x_1,\dots,x_d\}$, the one-dimensional locality is captured in the partitions $\mathcal{P}_k:=\{ \{x_1,\dots,x_k\} ,\{x_{k+1},\dots,x_d\} \}$, with respect to which one assumes low separation ranks $r_k:=r_{\mathcal{P}_k}$.
By reinterpreting the coefficient tensor $\theta$ as a matrix, via blocking the first $k$ tensor indices into a column index and the last $d-k$ tensor indices into a row index \cite{holtz_manifolds_2012}, and then applying the singular value decomposition (SVD), we find $r_k$ tensors $\alpha^{l_k}$ and $\beta^{l_k}$ such that
\begin{align}
\nonumber
    f^{\theta}&(x_1,\dots,x_d) :=\sum_{i_1,\dots, i_d} \theta_{i_1,\ldots, i_d}\psi_{i_1}(x_1)\cdot \dots \cdot \psi_{i_d}(x_d) \\
    & = \sum_{l_k=1}^{r_k} f^{\alpha_{l_k}}(x_1,\dots,x_k) \cdot f^{\beta_{l_k}}(x_{k+1},\dots, x_d) \, .
    \label{eq:alphabeta}
\end{align}
The SVD of $\theta$ therefore provides a decomposition of $f^{\theta}$ into a sum of separable functions. Furthermore, performing the SVD iteratively, for increasing $k$, results in a representation of ${\theta\in \mathbb{R}^{p}}$ via \emph{core tensors} ${A^k \in \mathbb{R}^{r_{k-1}\times \tilde{p} \times r_k}}$ \cite{oseledets_tensor-train_2011}.
The collection of core tensors $\{A^k\}_{k=1}^d$, contracted along their common indices indicated in Fig.~\ref{fig:tensortrain}c, is called a \emph{tensor train (TT) representation} of $\theta$.
The minimum cardinalities $r_k$ of the respective indices $l_k$ for a TT representation of $\theta$ are equal to the separation ranks in Eq. \eqref{eq:alphabeta}.
It has been further shown in \cite{holtz_manifolds_2012}, that the tensors $\theta$, which are expressible in the tensor train format for given ranks $(r_1,\dots,r_{d-1})$, build a sub-manifold $\Theta$ in the tensor space $\mathbb{R}^{\tilde{p}^d}$ with dimension
\begin{align}
\dim \Theta=\sum_{k=1}^d r_{k-1}\tilde{p} r_k -\sum_{k=1}^{d-1}r_k^2 \ll \tilde{p}^d = \dim \mathbb{R}^{\tilde{p}^d} \, ,
\end{align}
where we set ${r_0=r_d=1}$. This inequality holds, if $d$ is large and the ranks $r_k$ are bounded. 
Representing a tensor in the tensor train format with small ranks provides a way to store $\theta$ by its tensor train cores $A^k$. 
Assuming constant $r_k$, the memory demand scales linearly in $d$, which beats the `curse of dimensionality' manifested in the exponential scaling of the tensor space dimension.

A particular dynamical system featuring one-dimensional locality is described by the \emph{Fermi-Pasta-Ulam-Tsingou (FPUT)} 
equation, which models a vibrating string as a one-dimensional chain of masses. The variables of the system are the distortions of the masses, which are accelerated through nearest-neighbor interactions
\begin{align}
    \label{eq:FPU}
    \frac{d^2x_l(t)}{dt^2} =f_l(x_1& ,\dots,x_d)=  (x_{l+1}-2x_l+x_{l-1}) \\ &
    +\beta (x_{l+1}-x_l)^3 - \beta (x_{l}-x_{l-1})^3 \, .
    \nonumber
\end{align} 
As we will discuss in App.~\ref{app:fpu}, each tensor $\theta_l$, encoding the component $f_l$ of the FPUT equation, can be exactly decomposed as a TT with ranks $r_{l-1}=r_l=4$ and $r_{k}=1$ for $|k-l|>1$ , which is a direct consequence of the locality manifested in the nearest-neighbor interactions.  
This structural insight is also exploited in the work of \citet{gels_multidimensional_2019}.
The structure of the FPUT equation that leads to a restricted TT rank of $4$ is also present in a much larger class of functions, namely those that are quadratic forms of monomials of total degree smaller than $\tilde p$.
This condition ensures that the terms involving the largest degree monomials are univariate.
We summarize the statement for this class of functions as the following observation.

\begin{rem}[A function class with constant TT-rank]
\label{rem:lowrankpolynomials} 
Let $f: \mathbb R^{d} \to \mathbb R$ be a function of the form
\begin{equation*}
    f(x) = \sum_{l=1}^{d-1} \sum_{i+\tilde{i} < \tilde{p}} a(i,\tilde{i},l)\cdot  x_l^i x_{l+1}^{\tilde{i}}
\end{equation*}
where $a(i\textbf{},\tilde{i},l)\in\mathbb{R}$ are arbitrary coefficients. 
Further, let $f$ be represented by $\theta \in \mathbb R^{\tilde p \times \ldots \times \tilde p}$ using the dictionary consisting of products of the monomials $x_l^i, i = 0,\dots , \tilde{p}-1$.
Then, $\theta$ is of the form 
\begin{align*}
    \theta=\sum_{l=1}^{d-1}\sum_{i+\tilde{i} < \tilde{p}} a(i,\tilde{i},l) \cdot 1_{l-1}\otimes e^{(l)}_i\otimes e^{(l+1)}_{\tilde{i}}  \otimes 1_{d-l-1},
\end{align*}
in particular, the TT ranks of $\theta$ are bounded by $\tilde{p}$.
\end{rem}

The proof of this statement is a straight-forward calculation of the function form, Fig.~\ref{fig:tensortrain}, using the given expression for $\theta$. Similar results hold for other univariate dictionaries ${\psi_i :\mathbb{R}\to\mathbb{R}, i = 0,\dots,p-1}$, which include the constant function $1$.

\subsection{Correlations in the governing equation}\label{subsec:systems}

\begin{figure}
    \centering
  \begin{tikzpicture}[scale=0.3,thick]
\begin{scope}[shift={(20,0)},xscale=1,yscale=1]

    \coordinate[label=above:$c)$] (A) at (-3.5,0.5);

\draw (-1,-1) rectangle (1,1);
 \coordinate[label=above:$A^d$] (A) at (0,-1); 
 \draw (-1,0)--(-2.5,0) node[midway,above] {$i_d$};
  \draw (0,1)--(0,2.5) node[midway,right] {$l$};
 
  \draw (0,-1)--(0,-2) node[midway,right] {$l_{d-1}$};
   \draw (0,-3.5)--(0,-4.5) node[midway,right] {$l_2$};
    \draw (0,-6.5)--(0,-8) node[midway,right] {$l_1$};
  
\begin{scope}[shift={(0,-5.5)}]  
\draw (-1,-1) rectangle (1,1);
 \coordinate[label=above:$A^2$] (A) at (0,-1);
 \draw (-1,0)--(-2.5,0) node[midway,above] {$i_2$};
\end{scope}

  \coordinate[label=above:$\vdots$] (A) at (0,-3.7);
  
\begin{scope}[shift={(0,-9)}]  
\draw (-1,-1) rectangle (1,1);
 \coordinate[label=above:$A^1$] (A) at (0,-1);
 \draw (-1,0)--(-2.5,0) node[midway,above] {$i_1$};
\end{scope}

\end{scope}

    \coordinate[label=above:$a)$] (A) at (-3.5,0.5);
    
\draw (-1,-1) rectangle (1,1);
 \coordinate[label=above:$A^d$] (A) at (0,-1); 
 \draw (-1,0)--(-2.5,0) node[midway,above] {$i_d$};
 
  \draw (0,-1)--(0,-2) node[midway,right] {$l_{d-1}$};
   \draw (0,-3.5)--(0,-4.5) node[midway,right] {$l_2$};
    \draw (0,-6.5)--(0,-8) node[midway,right] {$l_1$};
  
\begin{scope}[shift={(0,-5.5)}]  
\draw (-1,-1) rectangle (1,1);
 \coordinate[label=above:$A^2$] (A) at (0,-1);
 \draw (-1,0)--(-2.5,0) node[midway,above] {$i_2$};
\end{scope}

  \coordinate[label=above:$\vdots$] (A) at (0,-3.7);
  
\begin{scope}[shift={(0,-9)}]  
\draw (-1,-1) rectangle (1,1);
 \coordinate[label=above:$A^1$] (A) at (0,-1);
 \draw (-1,0)--(-2.5,0) node[midway,above] {$i_1$};
\end{scope}

 \begin{scope}[shift={(-2.5,0)}]  
     \draw[fill=black] (7,-4.5) circle (0.25cm);
     \coordinate[label=left:$ \delta$] (A) at (7,-4);
      \coordinate[label=right:$l$] (A) at (7.5,-3.7);
    \draw (7,-4.5) to[bend right=70] (3.5,0);
    \draw (7,-4.5) to[bend left=20] (3.5,-5.5);
    \draw (7,-4.5) to[bend left=70] (3.5,-9);
    \draw (7,-4.5)--(9,-4.5);
 \end{scope}

 \begin{scope}[shift={(10,0)}]
 
     \coordinate[label=above:$b)$] (A) at (-3.5,0.5);
 
 \draw (-1,-1) rectangle (1,1);
 \coordinate[label=above:$A^d$] (A) at (0,-1); 
 \draw (-1,0)--(-2.5,0) node[midway,above] {$i_d$};
  \draw (1,0)--(2.5,0) node[midway,above] {$q_d$};
 
  \draw (0,-1)--(0,-2) node[midway,right] {$l_{d-1}$};
   \draw (0,-3.5)--(0,-4.5) node[midway,right] {$l_2$};
    \draw (0,-6.5)--(0,-8) node[midway,right] {$l_1$};
  
\begin{scope}[shift={(0,-5.5)}]  
\draw (-1,-1) rectangle (1,1);
 \coordinate[label=above:$A^2$] (A) at (0,-1);
 \draw (-1,0)--(-2.5,0) node[midway,above] {$i_2$};
 \draw (1,0)--(2.5,0) node[midway,above] {$q_2$};
\end{scope}

  \coordinate[label=above:$\vdots$] (A) at (0,-3.7);
  
\begin{scope}[shift={(0,-9)}]  
\draw (-1,-1) rectangle (1,1);
 \coordinate[label=above:$A^1$] (A) at (0,-1);
 \draw (-1,0)--(-2.5,0) node[midway,above] {$i_1$};
  \draw (1,0)--(2.5,0) node[midway,above] {$q_1$};
\end{scope}

\draw (2.5,1) rectangle (4.5,-10);
\coordinate[label=above:$S$] (A) at (3.5,-5.5);
\draw (4.5,-4.5) -- (6,-4.5)  node[midway,above] {$l$};
\end{scope}

\end{tikzpicture}
    \caption{Tensor network models for a system of functions $[f_1\dots f_d]$. In all models, for fixed equation index $l$, a one-dimensional locality assumption is made, such that $\theta_l$ can be represented in the TT format. (a) Under the assumption of independent equation components $\{f_l\}$, each tensor core inherits an additional index encoding the corresponding equation. These additional indices are contracted with a $\delta$ tensor which selects the appropriate TT for each equation component. (b) Under the assumption of a common structure between the TT for each equation component, a selection tensor $S$ can be used to exploit this structure by selecting the specific cores to be used for each component. (c) When a common structure is assumed, but not known beforehand, we add a single additional equation component index to the last TT core.}
    \label{fig:multiTT}
\end{figure}
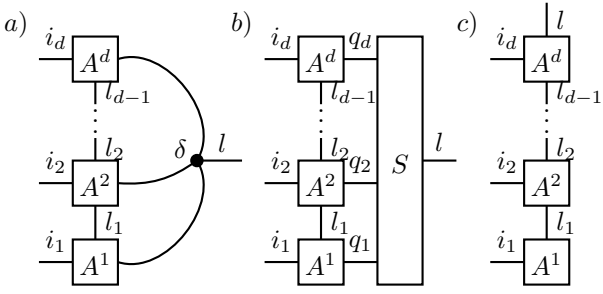

The recovery of a multi-component governing equation $[f_1\dots f_d]$ can be formulated as independently solving the tensor recovery problems  for each $\theta_l$. 
As per the previous section, the assumption of given separation
ranks for each tensor $\theta_l$ leads to an independent TT decomposition of each. We capture this structure in Fig.~\ref{fig:multiTT}a by a collection of $d$ different TT cores at each position, resulting in an additional index $l$ for each core tensor. 
By contraction with the delta tensor $\delta$, which selects for each equation $l$ the appropriate cores of the TT decomposition, the governing equation is represented by a single tensor $\theta \in \mathbb{R}^{\tilde{p} \times d}$.

Many governing equations however show correlations within the functions $f_l$ describing the dynamics of single variables, which are not taken into account by the model Fig.~\ref{fig:multiTT}a.
To illustrate this we recall the FPUT equation \eqref{eq:FPU}, in which each variable $x_k$ is represented in the system $[f_1 \dots f_d]$ in at most $n=4$ different ways.
We distinguish between the cases $k=l$, $k=l-1$, $k=l+1$ and $|k-l|>1$, which we refer to as the \emph{activation types} of the variable $x_k$. 
Given such correlation structure, there is a redundancy in the parameterization of the model by independent tensor trains for each function (Fig.~\ref{fig:multiTT}a), resulting from the appearance of identical\footnote{Up to symmetries in their parameterization resulting from transformations with invertible matrices \cite{holtz_manifolds_2012}.} 
TT cores representing different function tensors $\theta_l$ (see App.~\ref{app:fpu} for their explicit structure). Assuming the knowledge of the activation pattern, that is the activation type $i_{kl}\in \{1,\ldots,n\}$ of the variable $x_k$ in each function $f_l$, we introduce the model-specific selection tensor
\begin{align}
\label{def:FPU_Selection}
    S=\sum_{l=1}^d e^{(1)}_{i_{1l}}\otimes\cdots\otimes e^{(d)}_{i_{dl}} \otimes e_l \, \in \mathbb{R}^{n^d \times d } \,,
\end{align}
where we denote by $e^{(k)}_i$ the $i$'th basis vector in the $k$'th leg space of the tensor.
To remove the redundancies resulting from small numbers of activation types, we replace the function index $l$ at each tensor core $A^k$ by the activation type of each variable. If the number $n$ of activation types is small compared to the number $d$ of functions, i.e., if $d>4$ in the FPUT equation, this reparameterization results in an exponential decrease of the parameters.
The tensor $\theta$ representing the governing equation is then decomposed into the activation-informed tensor cores $A^k$ and the selection tensor $S$ (Fig.~\ref{fig:multiTT}b). 

In many settings, one might not have precise knowledge about the structure of the selection tensor, and can only assume the existence of a small number of activation types for each variable. In this situation one can parameterize the governing equation by a TT with an additional tensor leg representing the function index at one core tensor (Fig.~\ref{fig:multiTT}c). 
The thereby caused increase of the TT ranks is quantified in the next section for an exemplary class of systems.

\subsection{One-dimensional interacting systems}\label{sec:nn}

We now demonstrate the above parameterization methods for classes of governing equations, which describe systems with a one-dimensional locality structure in their degrees of freedom.

\begin{definition}[One-dimensional interacting systems] \label{def:nn}
We say that a governing equation $[f_1\dots f_d]$ is one-dimensional with \emph{interaction range} $(s_1,s_2)$ and \emph{separation rank} $N$, if there exist 
a set of functions $\{g_{i}(x)\}^{\bar p}_{i=1}$
and for each $l \in \{1, \ldots, d\}$ there is an index set $\mathcal{I}_l \subset \{1, \ldots, \bar p\}^{s_1 + s_2 + 1}$ with $|\mathcal{I}_l|\leq N$, such that
\begin{align*}
    f_l(x)=\sum_{(i_{l-s_1},\ldots, i_{l+s_2})\in \mathcal{I}_l} g_{i_{l-s_1}}(x_{l-s_1})\cdot \ldots \cdot g_{i_{l+s_2}}(x_{l+s_2})
\end{align*}
where we set $g_{i_k}=1$ for $k\leq 0$ and $k\geq d+1$. We further call the univariate function spaces $V$ spanned by the functions $g_{i}$ and including the constant function $1$ the \emph{leg embedding space} of $f$.
\end{definition}

For this classes of dynamical systems we derive the required ranks for the parameterization schemes, which have been discussed in Sec.~\ref{subsec:systems}.
\begin{theorem}[Tensor representation of one-dimensional interacting systems] \label{the:NN}
 Let $[f_1\dots f_d]$ be a one-dimensional governing equation with separation rank $N$ and interaction range $(s_1,s_2)$ and $V$ its leg embedding space.
Then, we have $f_l \in V^{\otimes d}$ for $l=1,\dots,d$ and:
\begin{enumerate}[(i)]
    \item The $\mathrm{CP}$-rank of each $f_l$ is at most $N$.
    \item Each $f_l$ has a TT representation, Fig.~\ref{fig:multiTT}a, with ranks $r_k\leq N$ for $l-s_1\leq k < l+s_2$ and $r_k=1$ else. 
    \item There is a selection tensor $S$ representation, Fig.~\ref{fig:multiTT}b, for $[f_1\ldots f_d]$ with $n=s_1+s_2+2$ and TT ranks bounded by $N$.
    \item The system $[f_1\dots f_d]$ has a single TT representation, Fig~\ref{fig:multiTT}c, with ranks $r_k\leq k-s_2+1+N(s_1+s_2)$.
\end{enumerate}
\end{theorem}
We give the proof of Thm.~\ref{the:NN} in App.~\ref{app:nn}.  
The theorem shows that each $f_l$ in a one-dimensional interacting system has bounded TT ranks independent of the number of variables $d$. 
Using a selection tensor $S$ to represent $[f_1\dots f_d]$ does not enlarge the TT ranks and allows for a total number of model parameters linear in $d$. 
Employing the single TT model instead of an explicit selection tensor, one finds a linear increase in the tensor train rank. As a consequence, the number of parameters scales as $d^3$ instead of $d$.

\section{Optimization of tensor networks}\label{sec:optimization}

We present in this section efficient learning algorithms for the proposed TT-type models. We begin by choosing an $\ell_2$-norm loss-function and utilizing Eq. \eqref{DefCoefficientRecovery}, which allows one to restate Problem~\ref{prob:GoverningEquation} as a least-squares optimization problem in the parameterization $\theta$:
\begin{equation}\label{eq:ptheta} \tag{$P_{\Theta}$}
   \operatorname*{minimize}_{\tilde \theta \in \mathbb R^{\tilde p \times \ldots \times \tilde p \times d}} \| 
    \Phi \tilde{\theta} - y \|_F^2 \quad \text{subject to}\quad  \tilde{\theta} \in  \Theta .
\end{equation}
Here, we restricted the optimization to the function space encoded by a subset of tensors $\Theta \subset \mathbb R^{\tilde p \times \ldots \times \tilde p \times d}$, which is taken to be one of the tensor network models discussed in the previous section, and illustrated in Fig.~\ref{fig:multiTT}.
For all of the models Problem \eqref{eq:ptheta} is an optimization task over $\{\tilde A^k\}_k$, the TT cores that determine $\tilde \theta$.
Hereinafter, we also use the symbol $f$ for the function defined by
$f(\tilde{A}^1,\ldots,\tilde{A}^d)\coloneqq \Phi\tilde{\theta}$, {i.e.}
the value of the candidate function for $f$ in the governing equation determined by the TT cores $\{\tilde A_k\}_k$ (evaluated at the given observations $x^j$).

\begin{figure}
    \centering
   \begin{tikzpicture}[scale=0.3,thick]

    \draw (-7.5,1) rectangle (-3.5,-10); 
    \coordinate[label=above:{$ \nabla_{A^2} f$}] (A) at (-5.5,-5.5);
   \draw (-3.5,-4.5)--(-2,-4.5) node[midway,above] {$j$};
   \draw (-7.5,-4.5)--(-9,-4.5) node[midway,above] {$i_2$};
   \draw (-7.5,-1.5)--(-9,-1.5) node[midway,above] {$l_2$};
      \draw (-7.5,-7.5)--(-9,-7.5) node[midway,above] {$l_{1}$};
       \draw[dashed] (-9,-1) rectangle (-11,-8); 
       \coordinate[label=above:{$A^2$}] (T) at (-10,-5.5);
      
\draw (-5.5,1)--(-5.5,2.5) node[midway,left] {$l$};
\coordinate[label=left:{=}] (T) at (0,-4.5);

\begin{scope}[yscale=-1,xscale=1,shift={(5,9)}]

    \draw (-1,-1) rectangle (3.5,1);
     \coordinate[label=below:$\psi_{i_1}(x^{\, j}_1)$] (A) at (1.25,-1.3);
     \draw (-1,0)--(-2.5,0) node[midway,above] {$i_1$};

    \begin{scope}[shift={(0,-3)}]  
       \draw (-1,-1) rectangle (3.5,1);
        \coordinate[label=below:$\psi_{i_2}(x^{\, j}_2)$] (A) at (1.25,-1.3);
        \draw (-1,0)--(-2.5,0) node[midway,above] {$i_2$};
    \end{scope}

    \coordinate[label=below:$\vdots$] (A) at (1.25,-7.5);
    \coordinate[label=below:$\vdots$] (A) at (5,-7.5);

    \begin{scope}[shift={(0,-9)}]  
        \draw (-1,-1) rectangle (3.5,1);
        \coordinate[label=below:$\psi_{i_d}(x^{\, j}_d)$] (A) at (1.25,-1.3);
        \draw (-1,0)--(-2.5,0) node[midway,above] {$i_d$};
    \end{scope}

    \draw[fill=black] (7,-4.5) circle (0.25cm);
     \coordinate[label=left:$ \delta$] (A) at (7,-5);
      \coordinate[label=right:$j$] (A) at (7.8,-5.5);
    \draw (7,-4.5) to[bend right=70] (3.5,0);
    \draw (7,-4.5) to[bend right=20] (3.5,-3);
    \draw (7,-4.5) to[bend left=70] (3.5,-9);
    \draw (7,-4.5)--(9,-4.5);

   \begin{scope}[shift={(-3.5,0)}]
\draw (-1,-1) rectangle (1,1);
 \coordinate[label=below:$\tilde{A}^1$] (A) at (0,-0.9);
  \draw (0,-1)--(0,-2) node[midway,left] {$l_{1}$};

    \draw (0,-6.75)--(0,-8) node[midway,left] {$l_d$};
  
\begin{scope}[shift={(0,-3)}]  
\draw[dashed] (-1,-1) rectangle (1,1);
 \coordinate[label=below:$A^2$] (A) at (0,-0.9);
   \draw (0,-1)--(0,-2.25) node[midway,left] {$l_2$};
\end{scope}

  \coordinate[label=above:$\vdots$] (A) at (0,-5);
\begin{scope}[shift={(0,-9)}]  
\draw (-1,-1) rectangle (1,1);
 \coordinate[label=below:$\tilde{A}^d$] (A) at (0,-0.9);
  \draw (0,-1)--(0,-2) node[midway,left] {$l$};  
\end{scope}
 \end{scope}
  \end{scope}
  
\end{tikzpicture}
    \caption{Gradient of the function $f(\tilde{A}^1,\ldots,\tilde{A}^d)$ with respect to a tensor train core 
    $A^2$, computed by a tensor network contraction.}
    \label{fig:gradient}
\end{figure}
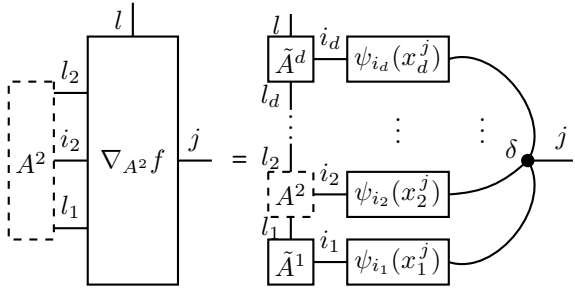

\subsection{Regularized alternating least squares}\label{subsec:ALS}

The \emph{alternating least squares (ALS)} strategy \cite{holtz_alternating_2012}
solves Problem~\eqref{eq:ptheta} by alternatingly optimizing the single TT cores while regarding the other TT cores as constant.
In each update step of the ALS algorithm, a tensor core $\tilde A^k$ is updated to a solution of
\begin{align}\label{eq:pk} \tag{$P_k$}
    \operatorname*{minimize}_{ A^k \in\mathbb{R}^{r_{k-1}\times \tilde{p} \times r_k}} \| f (\tilde{A}^1,\dots,A^k,\dots,\tilde{A}^d) - y \|_F^2.
\end{align}
After initialization by a random tensor in $\Theta$, the algorithm iterates so-called \emph{sweeps} in which the update step is performed for every tensor core, 
reminiscent of the density matrix renormalization group approach
\cite{DMRGWhite92}.
The optimal update can be found by the minimum criterion of first order, i.e., as the solution of the linear equation 
\begin{align}
\nonumber
    0&=\nabla_{A^k}\| f (\tilde{A}^1,\dots,A^k,\dots,\tilde{A}^d) - y \|_F^2 \, .
\end{align}
Due to the multi-linear dependence of the function $f$ on the TT cores and the measurement tensor $\Phi$, the analytical expression of the gradient $\nabla_{A^k}f$ is given by the tensor obtained via contraction of the measurement tensor $\Phi$ with all TT cores except for $A^k$ (see Fig.~\ref{fig:gradient}).
Crucially for the tensor networks considered here, and in fact for any tensor network with a tree structure, this contraction can be efficiently calculated.
As explained in App.~\ref{app:algorithm} each update step can thus be performed with computational complexity of $\mathcal{O}(r^3_{k-1}r^3_{k} {\tilde p}^3)$. 
Furthermore, the selection format introduces a non-uniqueness of the optimal solution for a fixed rank. We therefore regularize the ALS optimization problem by penalizing a large Frobenius norm of the TT core, see App.~\ref{subsec:normreg}.

\subsection{Rank-adaptive algorithm}\label{subsec:SALSA}
So far our formulation of the ALS algorithm relied on the separation ranks to be known. 
Physical principles such as locality might often merely justify the existence of a low-rank description, without providing precise upper bounds.
In order for the framework presented here to be broadly applicable, it is thus necessary to utilize an optimization algorithm which is able to adaptively identify admissible low-ranks through the course of optimization.
For the special case of the TT format, we demonstrate that the SALSA variant of the ALS algorithm by \citet{grasedyck_stable_2019} is able to provide such rank-adaptivity.
The major obstacle is that the solution of the ALS update step \eqref{eq:pk} is unstable under changes of the TT-rank. To resolve this instability problem, given a fixed $k$, \citet{grasedyck_stable_2019} introduce a unique decomposition of a tensor $\theta$ into orthogonal tensors $\mathcal{L},\mathcal{R}$ and diagonal matrices $\Sigma_{\mathcal{L}},\Sigma_{\mathcal{R}}$,
\begin{center}
\vspace{.3ex}
   \begin{tikzpicture}[scale=0.3]

\draw (0,0) rectangle (5,2);
 \coordinate[label=above:$ \theta$] (A) at (2.5,0.15);
\draw (0.5,0)--(0.5,-1.2) node[midway,left] {$i_{\mathcal{L}}$};
\draw (2.5,0)--(2.5,-1.2) node[midway,left] {$i_k$};
\draw (4.5,0)--(4.5,-1.2) node[midway,left] {$i_{\mathcal{R}}$};

 \coordinate[label=above:{$=$}] (A) at (6,0.2);
\draw (7,0) rectangle (9,2);
\coordinate[label=above:$ \mathcal{L}$] (A) at (8,0.15);
\draw (8,0)--(8,-1.2) node[midway,left] {$i_{\mathcal{L}}$};
\draw (9,1)--(10,1); 

\draw (10,0) rectangle (12,2);
\coordinate[label=above:$\Sigma_{\mathcal{L}}$] (A) at (11,0.05);
\draw (12,1)--(13,1); 

\draw (13,0) rectangle (15,2);
\coordinate[label=above:$ \mathcal{N}^k$] (A) at (14,0.1);
\draw (14,0)--(14,-1.2) node[midway,left] {$i_k$};
\draw (15,1)--(16,1); 

\draw (16,0) rectangle (18,2);
\coordinate[label=above:$ \Sigma_{\mathcal{R}}$] (A) at (17,0.05);
\draw (18,1)--(19,1); 

\draw (19,0) rectangle (21,2);
\coordinate[label=above:$ \mathcal{R}$] (A) at (20,0.15);
\draw (20,0)--(20,-1.2) node[midway,left] {$i_{\mathcal{R}}$};

\end{tikzpicture} ,
   \vspace{.2ex}
\end{center}
such that $\mathcal{L} \Sigma_{\mathcal{L}} \big[ \mathcal{N} \Sigma_{\mathcal{R}}  \mathcal{R} \big] $ and $\big[ \mathcal{L} \Sigma_{\mathcal{L}} \mathcal{N} \big] \Sigma_{\mathcal{R}}  \mathcal{R} $ are singular value decompositions. 
To stabilize the update  
\eqref{eq:pk}, one replaces the objective function by its average over a local neighborhood of $\theta$ with diameter $\omega$.
This is equivalent to the addition of regularization terms to the optimization problem:
\begin{align}\label{eq:salsa}\tag{$P_{k-\mathrm{s}}$}
    \operatorname*{minimize}_{\mathcal{N}^k\in\mathbb{R}^{r_{k-1}\times \tilde{p} \times r_k}} 
        &\| f(\mathcal{L},\mathcal{N}^k,\mathcal{R}) - y \|_F^2 \\
    &\quad+ \omega^2\left(\|\Sigma_{\mathcal{L},\epsilon}^{-1}\mathcal{N}^k\|_F^2 + \|\mathcal{N}^k\Sigma_{\mathcal{R},\epsilon}^{-1}\|_F^2\right)
    . \nonumber
\end{align}
For the numerical inversion of $\Sigma_{\mathcal{L}}$ and $\Sigma_{\mathcal{R}}$, the singular values below a certain threshold $\epsilon$ are set to $\epsilon$. 
The rank-adaption strategy enforces a constant number $r_{\mathrm{min}}$ of singular values to be below $\epsilon$ for $\Sigma_{\mathcal{L}}$ and $\Sigma_{\mathcal{R}}$ in the following way: If after solving \eqref{eq:salsa} the number of singular values below $\epsilon$ increases or decreases, a corresponding number of randomly chosen singular vectors in $\mathcal{L}$ or in $\mathcal{R}$ are discarded or
added, respectively. 
The resulting SALSA algorithm as well as a detailed description of the chosen hyper-parameters is given in the Appendix.

\section{Numerical simulations} \label{sec:examples}

Finally, we demonstrate the performance of ALS and its rank-adaptive generalization SALSA for the introduced tensor network parameterizations.
To this end, we study variants of the Fermi–Pasta–Ulam–Tsingou (FPUT) equation, as introduced in Eq.~\eqref{eq:FPU}, using a univariate dictionary consisting of the first ${\tilde{p}=4}$ $L_2$-orthogonal Legendre polynomials. 

We compare the recovery of the coefficient $\theta$ for two models: One involving the selection tensor $S$ (Fig.~\ref{fig:multiTT}b) and the other regarding the single TT representation (Fig.~\ref{fig:multiTT}c). 
We use the regularized ALS algorithm to recover the model involving the selection tensor.  The selection tensor has unit-rank tensors in its null space preventing rank-adaption in the sense of SALSA.
Using the rank-adaptive SALSA, the single TT model is optimized without specifying the TT ranks.

\begin{figure}[tb]
    \vspace{-1.2cm}
    \centering
    \input{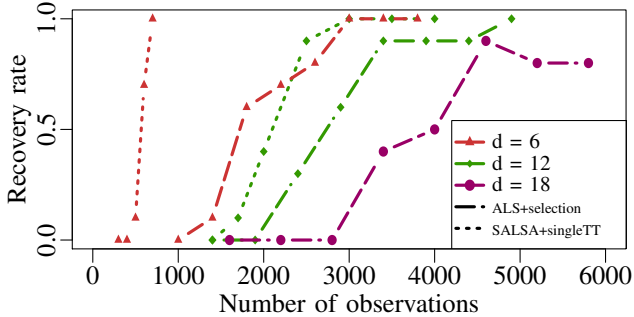}
    \vspace{-1cm}
    \caption{\label{fig:stFPUT}
    The recovery rate (relative error threshold $10^{-6}$) for different number of observations $m$ and dimension $d$ for the FPUT equation with $\beta = 0.7$ using the selection tensor model (Fig.~\ref{fig:multiTT}b), trained with ALS for $15$ iterations  (long dash lines), and the single TT model (Fig.~\ref{fig:multiTT}c), trained with SALSA for $60$ iterations (dotted lines).
    Each point is the average value of $10$ trials.
    }
    \vspace{-.3cm}
\end{figure}

We consider a coefficient tensor as successfully recovered if it relatively deviates from the correct coefficient tensor by less then $10^{-6}$ in Frobenius norm.
The respective recovery rates for both models over ten experiments, for different values of $d$ and $m$, are compared in Fig.~\ref{fig:stFPUT}.
As expected, the required number of observations in both parameterization formats increases with the number of variables $d$.   
As per Thm.~\ref{the:NN}, note that while the selection tensor $S$ allows for a model with constant TT ranks, the single TT model requires TT ranks that linearly increase with the position of the core tensors. 
Details concerning the exact ranks are provided in App.~\ref{app:fpu}.
For small $d$, instances of the single TT model are already successfully recovered for smaller observation numbers compared to the selection tensor model, which reflects the more efficient parameterization with a single TT in this regime. 
The SALSA algorithm furthermore recovers the exact ranks of the model.
For high $d$, however, increasing ranks of the single TT model cause an inadmissible computational demand, which scales for each update step quadratically in the dimension of the involved core tensor. 
Correspondingly, we observe in our numerics that performing the SALSA method on the single TT model becomes intractable on desktop hardware. 
In this regime, the model involving a selection tensor, with bounded ranks, is less computationally demanding. 
In conclusion, we find a trade-off between the computational cost and the amount of prior knowledge about the equations structure that the different models encode.

An extension of the FPUT equation that is still representable by the selection tensor model (Fig.~\ref{fig:multiTT}b) includes additional non-local mean field terms.
ALS can be successfully employed for the recovery of this type of equations. We give detailed numerical results in the App.~\ref{app:numdiscuss}. 
As the FPUT model is rather restrictive, we also studied the recovery of arbitrary one-dimensional interacting systems with interaction range $(1,1)$ and random coefficients. 
As detailed in App.~\ref{app:numdiscuss},
we still achieve adequate recovery rates with moderately sized numbers of observations for the selection tensor model (Fig.~\ref{fig:multiTT}b) with $d=6$ and $d=12$. For $d=18$ the recovery becomes more involved due to local minima. 

All simulations were implemented in python using \texttt{xerus} \cite{xerus} and 
run on a computer with $4\times\SI{3500} {\mega\hertz}$ processors and $\SI{16}{\giga\byte}$ RAM. As an example for the selection tensor format a single simulation for $d=18$ and $m=6000$ takes about $\SI{500}{\second}$ for $20$ iterations. 

\section{Conclusion and outlook}
In this work, we have provided a tensor network based framework for learning non-linear dynamical laws for systems with many degrees of freedom, in which the incorporation of fundamental physical principles, such as locality of admissible correlations, leads naturally to an efficient and scalable approach.
In particular, we have shown that multivariate function dictionaries, built by products of univariate functions, give rise to tensor network structured representations for governing equations.
Furthermore, we have argued that low separation ranks of these tensor networks are in accordance with constraints on the correlation structure of a system's variables, and as such have provided a generic scheme to exploit expert knowledge on the structure of governing equations. 
As a guiding example, we have discussed multiple tensor network models for systems with a constrained one-dimensional correlation structure, 
and obtained model-specific rank bounds for a large class of physical systems.
We have furthermore adapted and implemented fixed-rank and rank-adaptive optimization schemes that does not require an a priori knowledge of the ranks. 
Our algorithms have numerically been demonstrated to be able to successfully recover variants of Fermi-Pasta-Ulam-Tsingou equation and random locally interacting systems. 
We regard our work as a fruitful step towards understanding the potential for tensor network parameterizations of governing equations in the light of fundamental physical principles.

In the context of many-body physics, the relation between constrained correlations fulfilling an \emph{area law} \cite{eisert_area_2010} and low-rank TT representations have been rigorously established \cite{Schuch_MPS}.
While we here only exemplify a similar relation for specific models, 
we expect rigorous information theoretical analyses of the correlation structure to yield analogous results, including notions of expressivities of
tensor networks capturing certain data structures \cite{Expressive}.
Additionally, there exist a variety of other tensor network model that can be used in order to extend our approach, such as projectively entangled pair states for the description of higher dimensional systems \cite{VerstraeteBig} or the 
multi-scale renormalization ansatz 
\cite{Vidal:2007:mera} that are expected to describe scale-invariant correlation structures.
Our approach, in principle, admits the implementation of symmetry concepts, by further constraining the structure of decomposing core tensors \cite{cohen_expressive_2015,Symmetries}. 
Additionally, the expressivity of our parameterization schemes is determined by univariate dictionary functions, and incorporating the selection of a suitable dictionary in the learning task along the ideas of \citet{champion_data-driven_2019-1} and \citet{PhysRevLett.124.010508} is another interesting extension.
From a compressed sensing perspective, rigorous tensor recovery guarantees have been proven  \cite{rauhut_tensor_2015-1,rauhut_low_2016,grotheer_iterative_2019}. 
It is an open problem to extend these results to the approach taken in this work. 

Finally, we would like to highlight that there are different perspectives on the recovery of dynamical systems that do not directly aim at the recovery of governing equations \cite{williams_datadriven_2015,kevrekidis_kernel-based_2016,kutz_dynamic_2016,klus_data-driven_2018}. 
It is an interesting question whether one can also use tensor network models that encode physical principles to improve the scalability of these alternative approaches
(see, for example, Ref.~\cite{nuske_tensor-based_2019} for recent progress in this direction).
It is even conceivable to include tensor structured coefficients to increase the expressivity and scalability of symbolic regression algorithms. 
It is the hope that the present work stimulates such endeavours.

\section*{Acknowledgements}
We are grateful to Patrick Gel{\ss}, Stefan Klus and Christof Sch{\"u}tte for many fruitful discussions and extensive explanations about the SINDy and MANDy approaches. 
A.~G. and G.~K.~acknowledge funding from the MATH+ project EF1-4.
M.~G.~acknowledges funding from the DFG (SCHN 530/15-1).
R.~S.~acknowledges funding from the Alexander von Humboldt foundation. 
I.~R.~and J.~E.~acknowledge funding from the DFG (EI 519/9-1 CoSIP, CRC 1114 project B06,
MATH+ project EF1-7, CRC 183 project B01, and EI 519/15-1) and 
the BMWi (PlanQK). 

\bibliographystyle{icml2020}



\onecolumn 
\appendices

\section*{\Large Appendices}
\vspace{1\baselineskip}

\noindent In these appendices we provide additional details concerning tensor network notation and the tensor representation of multivariate functions (App.~\ref{app:notation}), as well as a proof of the expressivity result Theorem ~\ref{the:NN} (App.~\ref{app:nn}). In  App.~\ref{app:govEqs} we further  provide details on the types of governing equations that were studied in this work. Then, in App.~\ref{app:numdiscuss} we provide additional numerical results, along with an extended discussion. Finally, in App.~\ref{app:algorithm} we provide a detailed account of the algorithmic implementation of all numerical methods. 

\section{Tensor notation}\label{app:notation}

We briefly describe here the graphical notation used in this work to represent both tensors and tensor network contractions.
This graphical notation is particularly common in the many-body physics literature \cite{VerstraeteBig,eisert_area_2010,Orus-AnnPhys-2014}, and we refer there for a more detailed presentation.
In general, an order $d$ tensor $\theta \in \mathbb{R}^{\tilde{p}^d}$ is represented by a block with $d$ different legs, each of which represents an index of the tensor,
\begin{equation}
\theta = \begin{tikzpicture}[scale=0.3,baseline=-3.5pt]]
    \draw (-8,-1) rectangle (-1,1);
    \coordinate[label=above:{$\theta$}] (A) at (-4.5,-0.9);
    \draw (-7.25,-1)--(-7.25,-2.25);
    \draw (-5.75,-1)--(-5.75,-2.25) ;
    \draw (-1.75,-1)--(-1.75,-2.25);
    \coordinate[label=below:{$\dots$}] (A) at (-3.5,-1.2);
\end{tikzpicture} \, .
\end{equation}
The particular scalar elements $\theta_{i_1 \ldots i_d} \in \mathbb{R}$ of the tensor are then represented by indicating the appropriate indices at the relevant legs of the tensor, i.e.
\begin{equation}
\theta_{i_1\ldots i_d} = \begin{tikzpicture}[scale=0.3, baseline = -3.5pt]
    \draw (-8,-1) rectangle (-1,1);
    \coordinate[label=above:{$\theta$}] (A) at (-4.5,-0.9);
    \draw (-7.25,-1)--(-7.25,-2.25) node[midway,left] {$i_1$};
    \draw (-5.75,-1)--(-5.75,-2.25) node[midway,left] {$i_2$};
    \draw (-1.75,-1)--(-1.75,-2.25) node[midway,left] {$i_d$};
    \coordinate[label=below:{$\dots$}] (A) at (-4.25,-1.2);
\end{tikzpicture} \, .
\end{equation}
For example, a vector $\theta \in \mathbb{R}^p$ and a matrix $\Phi \in \mathbb{R}^{m\times p}$ would be represented by a box with one leg and two legs respectively,
\begin{equation}
    \theta_i = \begin{tikzpicture}[scale=0.3,baseline=-3.5pt]
    \draw (6,0)--(7.5,0) node[midway,above] {$i$};
    \draw (7.5,-1) rectangle (9.5,1);
    \coordinate[label=above:$\theta$] (A) at (8.5,-0.8);
\end{tikzpicture}\, \qquad \Phi_{ji} = \begin{tikzpicture}[scale=0.3,baseline=-3.5pt]
    \draw (6,0)--(7.5,0) node[midway,above] {$j$};
    \draw (7.5,-1) rectangle (9.5,1);
    \draw (9.5,0)--(11,0) node[midway,above] {$i$};
    \coordinate[label=above:$\Phi$] (A) at (8.5,-0.8);
\end{tikzpicture} \, .
\end{equation}
A contraction over a particular tensor index is then represented by connecting the appropriate legs of the relevant tensors. For example, matrix vector multiplication is indicated via

\begin{equation}
    (\Phi\theta)_j = \sum_{i}\Phi_{ji}\theta_i = \begin{tikzpicture}[scale=0.3,baseline=-3.5pt]
    \draw (2.5,0)--(4,0) node[midway,above] {$j$};
    \draw (4,-1) rectangle (6,1);
    \coordinate[label=above:$\Phi$] (A) at (5,-0.8);
    \draw (6,0)--(7.5,0) node[midway,above] {$i$};
    \draw (7.5,-1) rectangle (9.5,1);
    \coordinate[label=above:$\theta$] (A) at (8.5,-0.8);
\end{tikzpicture}\, .
\end{equation}
Occasionally, when it is clear which index is being contracted, we will omit the summation index above the corresponding connected legs.
This graphical notation therefore allows us to represent in a concise way tensors built from the contraction of component tensors, through networks of tensor diagrams connected in the appropriate fashion.

In this work we are primarily concerned with using tensor networks in order to represent multivariate functions, living in a linear space built from the tensor product of univariate function spaces. From a more general perspective consider a collection of linear spaces $\{V_k\}|_{k=1}^d$, corresponding with the univariate function spaces, each with a basis $\{e^{(k)}_i\}_{i=1}^p\subset V_k$. We can then define the tensor space $V = V_1\otimes\cdots\otimes V_d$ by its basis, given by the tensor products of basis elements for the component linear spaces, i.e.
\begin{align}
    e_{i_1,\ldots,i_d} =e^{(1)}_{i_1}\otimes \cdots \otimes e^{(d)}_{i_d} \, .
\end{align}
With respect to this basis, an arbitrary element $v \in V$ of the tensor product space can then be expressed by the coefficient tensor $\theta \in \mathbb{R}^{p^d}$ which stores the coefficients of the components of $v$, i.e.
\begin{equation}
    v = \sum_{i_1,\ldots,i_d}\theta_{i_1,\ldots,i_d}e_{i_1,\ldots,i_d} \, .
\end{equation}

To build multivariate function spaces, we consider as building blocks vector spaces $V_k$ of univariate functions, specified by basis functions $e^{(k)}_{i}\coloneqq \psi^{(k)}_{i}: \mathbb R \to \mathbb R, x \mapsto \psi^{(k)}_{i}(x)$. Multivariate functions are then constructed by tensor product of basis functions via
\begin{equation}
\label{DefTensorProduct}
    \psi^{(1)}_{i_1} \otimes \dots \otimes  \psi^{(d)}_{i_d}: (x_1, \dots , x_d) \mapsto \psi^{(1)}_{i_1}(x_1)\cdot \dots  \cdot \psi^{(d)}_{i_d}(x_d) \, .
\end{equation}
The tensor space $V = V_1\otimes\cdots\otimes V_d$ is then defined by the linear hull of these tensor products.
As a result, given $d$ univariate function spaces $\{V_k\}$, we can specify an arbitrary element $f:\mathbb{R}^d \rightarrow \mathbb{R}$ of the tensor space $V = V_1\otimes\cdots\otimes V_d$ via its coefficient tensor $\theta \in \mathbb{R}^{\tilde{p}^d}$. Explicitly, we have that
\begin{align}
    f(x_1,\ldots,x_d) &= \sum_{i_1,\ldots,i_d}\theta_{i_1,\ldots,i_d}\psi^{(1)}_{i_1}(x_1) \cdot \psi^{(2)}_{i_2}(x_2)\ldots \psi^{(d)}_{i_d}(x_d).
\end{align}
Finally, note that we employ the symbol $1_d$ to represent the constant function with unit value, within the space of multivariate functions acting on $d$ variables. It can be  constructed by the tensor products of the constant functions $1\in V_k$.

\section{Proof of Theorem \ref{the:NN}} \label{app:nn}
We provide here a proof of Theorem~\ref{the:NN}, which states rank bounds for the governing equation of a one-dimensional interacting system with range $(s_1, s_2)$ and separation rank $N$.
According to Def.~\ref{def:nn}, such a governing equation $[f_1 \ldots f_d]$ is of the form 
\begin{equation}
    \label{eq:NN}
    f_l(x)=\sum_{(i_{l-s_1},\ldots, i_{l+s_2})\in \mathcal{I}_l} g_{i_{l-s_1}}(x_{l-s_1})\cdot \ldots \cdot g_{i_{l+s_2}}(x_{l+s_2}).
\end{equation}
Let $V$ be the embedding space for $f_l$, i.e. an appropriate univariate function space from which to build the multivariate function space of which $f_l$ is an element. Specifically, as explained in App. ~\ref{app:notation}, we regard each function $g_{i_k}$ as an element of the vector space $V$ and thus $g_{i_1}\otimes \ldots \otimes g_{i_d}$ as a separable tensor in $V^{\otimes d}$.
In order to be more precise we can make the dependence of $f_l$ on all variables explicit by using the constant function (as described in App.~\ref{app:notation}) for all variables outside of the interaction range. Eq.~\ref{eq:NN} then corresponds to
\begin{align}
\label{eq:NNabstract}
    f_l =\sum_{(i_{l-s_1}\ldots i_{l+s_2})\in \mathcal{I}_l} 1_{l-s_1-1}\otimes g_{i_{l-s_1}}\otimes
    \cdots\otimes g_{i_{l+s_2}} \otimes 1_{d-l-s_2}\,.
\end{align}
This is a $\mathrm{CP}$-decomposition of the tensor $f_l$. 
Since we have $|\mathcal{I}_l|\leq N$, the $\mathrm{CP}$-rank of $f_l$ is therefore bounded by $N$, establishing the claim of (i).

For (ii) note that the minimal TT ranks for the representation of a tensor are equal to the separation ranks with respect to the partitions $\mathcal{P}_k=\{\{x_1,\ldots,x_k\}, \{x_{k+1},\ldots,x_d\} \}$.
A formal proof of this statement is given in  Ref.~\cite{holtz_manifolds_2012}.
If $k<l-s_1$ the decomposition \eqref{eq:NNabstract} is equal to 
\begin{align*}
    f_l= 1_{k}\otimes\left[ \sum_{(i_{l-s_1}, \ldots, i_{l+s_2})  \in \mathcal{I}_l} 1_{l-s_1-k-1} \otimes g_{i_{l-s_1}}\otimes\cdots \otimes g_{i_{l+s_2}} \otimes 1_{d-{l}-s_2} \right] \, ,
\end{align*}
thus the respective separation rank and, correspondingly, the TT rank $r_k$ is at most one. This bound holds also for $k\geq l+s_2$, which follows from the analogous decomposition 
\begin{align*}
    f_l=  \left[ \sum_{(i_{l-s_1}, \ldots, i_{l+s_2})\in \mathcal{I}_l} 1_{l-s_1-1} \otimes g_{i_{l-s_1}}\otimes \cdots \otimes g_{i_{l-s_2}} \otimes 1_{k-l-s_2} \right] \otimes 1_{d-k} \, .
\end{align*}
For $k \in \{ l-s_1,\ldots ,l+s_2-1\}$ $f_l$ does not permit such a factor decomposition in general
but remains of the form
\begin{align}
\nonumber
    f_l=\sum_{(i_{l-s_1},\ldots, i_{l+s_2})\in \mathcal{I}_l} \big[1_{l-s_1-1} \otimes g_{i_{l-s_1}}\otimes\cdots\otimes g_{i_{k}}\big] \otimes \big[g_{i_{k+1}}\otimes \cdots \otimes g_{i_{l+s_2}} \otimes 1_{d-{l}-s_2}\big] \, .
\end{align}
We thus have a separation rank bounded as $r_k\leq |\mathcal{I}_l|\leq N$.

For (iii) we observe that at most only
$s_1+s_2+1$ different functions $f_l$ depend non-trivially on the single variable $x_k$ for all $k$.
Thus, each variable $x_k$ has  $s_1+s_2+1$ different activation types.
The trivial dependence on a variable in all other functions constitutes another activation type. 
Hence, we need at most $n=s_1+s_2+2$ activation types to represent the variables for the system $[f_1 \dots f_d]$. 
We can thus build a selection tensor $S$ by
\begin{align}
\nonumber
    S=\sum_{l=1}^{d}  e^{(1)}_{1}\otimes \cdots \otimes e^{(l-s_1-1)}_{1} \otimes e^{(l-s_1)}_2 \otimes \cdots \otimes e^{(l+s_2)}_{s_1+s_2+2}  \otimes e^{(l+s_2+1)}_{1}\otimes \cdots \otimes e^{(d)}_{1} \otimes e_l \, .
\end{align}
Implementing an additional leg with dimension $n=s_1+s_2+2$ indexing the different active TT cores $A^k$ for each variable $x_k$ results in a representation of $[f_1\ldots f_d]$ by a contraction with $S$ (Fig.~\ref{fig:multiTT}b). At each position $k$ the required TT rank $r_k$ is therefore bounded by the maximum separation rank of the single functions $f_l$, which is by (ii) bounded by $N$.

(iv) We represent the governing equation $[f_1 \ldots f_d]$ by the tensor $f \in V^{\otimes d} \otimes \mathbb{R}^d$, 
\begin{align}
    \label{NNterm0}
    f   =&  \sum_{l=1}^d 1_{l-s_1-1} \otimes \tilde{f}_l \otimes 1_{d-l-s_2} \otimes e_l \, ,
\end{align}
where we denote by $\tilde{f}_l$ the projection of $f_l$ to the product space of its legs $l-s_1,\ldots,l+s_2$. 
Analogous to the proof of (ii), it is enough to bound the separation ranks $r_k$ of the tensor $f$ with respect to the partitions ${\mathcal{P}_k}$. For each $k=1,\ldots ,d$ we therefore split the sum \eqref{NNterm0} into the terms
\begin{align}
    \label{NNterm1}
f= &  \sum_{l \leq k-s_2} 1_{l-s_1-1} \otimes \tilde{f}_l  \otimes 1_{d-l-s_2} \otimes e_l\\
     \label{NNterm2}
    &+ 1_{k} \otimes \big[\sum_{l> k+s_1} 1_{l-k-s_1-1} \otimes \tilde{f}_l \otimes 1_{d-l-s_2} \otimes e_l \big] \\
    &+ \sum_{l=k-s_2+1}^{k+s_1} \sum_{(i_{l-s_1},\ldots, i_{l+s_2})\mathcal{I}_l} 1_{l-k-s_1-1} \otimes g_{i_{l-s_1}}\otimes\cdots \otimes g_{i_{l+s_2}} \otimes 1_{d-l-s_2} \otimes e_l \, .
        \label{NNterm3}
\end{align}
The tensor in term \eqref{NNterm1} is nonzero if $k-s_2>0$ and has a separation rank with respect to $\mathcal{P}_k$ of at most $k-s_2$. Term \eqref{NNterm2} is nonzero in case $k<d-s_1$ and then contributes with a separation rank of at most $1$. For each summand of the first sum of term \eqref{NNterm3} the corresponding tensor has a separation rank bounded by $N$, since $|\mathcal{I}_l|\leq N$. The number of these summands is given by $s_1+s_2-\max(s_1-k+1,0)-\max(k-d-s_2,0)$, thus bounded by $s_1+s_2$. In all cases we have thus observed a separation rank bound of $r_k \leq k-s_2+1+N(s_1+s_2)$. \hfill\qed

\section{Examples of governing equations} \label{app:govEqs}
We now discuss examples of governing equations, which illustrate the expressivity results and provide test systems for the numerical experiments. We follow closely the notation and techniques of \citet{gels_multidimensional_2019}.

\subsection{Variants of the Fermi-Pasta-Ulam-Tsingou equation} \label{app:fpu}

Firstly we discuss variants of the \emph{Fermi–Pasta–Ulam–Tsingou (FPUT)} equation, originally introduced by \citet{fermi_studies_1955}. Given constants $m_l$ and $\beta_l$ we define for $l = 1,\ldots, d$ the functions
\begin{align}\label{eq:FPUMean}
\begin{split}
    \frac{d^2}{dt^2} x_l(t) =& f_l(x(t))= \big( x_{l+1}(t) - 2x_l(t) + x_{l-1}(t) \big) + \beta_l\big(x_{l+1}(t) - x_l(t)\big)^3 - \beta_l\big( x_l(t)-x_{l-1}(t) \big)^3 + \sum_{\tilde{l}=1}^dm_{\tilde{l}}x_{\tilde{l}} \, .
\end{split}
\end{align}
where ${x_0 = x_{d+1} = 0}$. 
Note that enforcing $\beta_l = \beta$ for all $l$ results in a translationally invariant system of equations, and that this translation invariance can be broken by allowing for different values of $\beta_l$ for each $l$.
We will show in the following, that this does not affect the required ranks in the representation.
Furthermore, in addition to the interaction terms, we have included the constant field term $\sum_{\tilde{l}=1}^dm_{\tilde{l}}x_{\tilde{l}}$.

We now derive an explicit parameterization of the functions $f_l$ in the TT format (see Fig.~\ref{fig:tensortrain}c). While in our numerical studies we represented the equations \eqref{eq:FPUMean} with respect to $L_2([-1,1])$-orthonormal basis functions (App.~\ref{app:L2orth}), here we choose the monomials $\{1,x_k,x_k^2,x_k^3\}$ as basis functions $\{\psi_{i_k}\}$ for each variable $x_k$. As such, we remark that the representability of the coefficient tensor in the following tensor network formats is invariant under orthonormalization of the dictionary $\{\psi_{i}\}_{i=1}^{\tilde{p}}$ with respect to any scalar product.

The variable dependencies of Eq.~\eqref{eq:FPUMean} unravel the structure of the corresponding tensor $\theta_l$ to be
\begin{align}
\label{eq:FPUthetal}
    \theta_l=1_{l-2}\otimes \tilde{\theta}_l \otimes 1_{d-l-1} + \sum_{\tilde{l}=1}^d 1_{\tilde{l}-1}\otimes m_{\tilde{l}} e^{(\tilde{l})}_1 \otimes 1_{d-\tilde{l}}\, ,
\end{align}
where $\tilde{\theta}_l$ is a tensor of order three and again we have denoted by $1$ the constant function in the respective function spaces (see App.~\ref{app:notation}). 
Factorization of the terms in Eq.~\eqref{eq:FPUMean} corresponds to a $\mathrm{CP}$-decomposition of $\tilde{\theta}_l$ as
\begin{align*}
    \tilde{\theta}_l &= e^{(l-1)}_1 \otimes  \left[ [-2e^{(l)}_2-2\beta_l e^{(l)}_4] \otimes e^{(l+1)}_1 + [e^{(l)}_1+3 \beta_l e^{(l)}_3] \otimes e^{(l+1)}_2   - 3 \beta_l e^{(l)}_2 \otimes e^{(l+1)}_3+  \beta_l e^{(l)}_1 \otimes e^{(l+1)}_4 \right], \\
    &+ e^{(l-1)}_2 \otimes [ e^{(l)}_1+3 \beta_l e^{(l)}_3] \otimes e^{(l+1)}_1,\\
    &-e^{(l-1)}_3 \otimes 3 \beta_l e^{(l)}_2 \otimes e^{(l+1)}_1\\
    & + e^{(l-1)}_4 \otimes \beta_l e^{(l)}_1 \otimes e^{(l+1)}_1 \, 
    .
\end{align*}
To derive from the above $\mathrm{CP}$-decomposition an illustration of the TT format of $\tilde{\theta}_l$ we build a matrix $A^l$ by the coefficient vectors in the second leg space, which corresponds to the univariate function space of the variable $x_l$. Following Ref.~\cite{gels_multidimensional_2019} the matrix is
\begin{align}
\nonumber
A^l=&\left[ \begin{array}{rrrr}
-2e_2-2\beta_l e_4 & e_1 + \beta_l e_3 & -3\beta_l e_2 & \beta_l e_1 \\
e_1+3\beta_l e_3 & 0 & 0 & 0 \\
-3\beta_l e_2 & 0 & 0 & 0 \\
\beta_l e_1 & 0 & 0 & 0 \\
\end{array}\right]^{(l)}
\nonumber .
\end{align}
Defining further matrices ${A^{l-1}=[e_1^{(l-1)},e_2^{(l-1)},e_3^{(l-1)},e_4^{(l-1)}]}$ and ${A^{l+1}=[e_1^{(l+1)},e_2^{(l-1)},e_3^{(l-1)},e_4^{(l+1)}]^T}$, the tensor $\tilde{\theta}_l$ is given by the matrix contraction $A^{l-1}\cdot A^l \cdot A^{l+1}$, performed coordinate-wise as a tensor product (see Eq.~\eqref{DefTensorProduct}). 
We can now exploit the above tensor cores to build a TT decomposition of the full tensor $\theta_l$. Specifically, we see that this requires TT ranks of $r_k=4$ if $k=l-1,l$ and $r_k=2$ for all other $k$:
\begin{align}
\nonumber
\small
&\theta_{l}=  \left[ \begin{array}{rr}e_1&m_1e_2\end{array}\right]^{(1)} \cdots \left[\begin{array}{rr}e1&m_{l-2}e_2\\0&e_1\end{array}\right]^{(l-2)}  \cdot  \\
& \nonumber \small \left[ \begin{array}{rrrr}
e_1 & e_2 & e_3 & e_4 \\ 0 & 0 & 0 & \beta_l^{-1}e_1
\end{array}\right]^{(l-1)}\cdot 
\left[ \begin{array}{rrrr}
{(-2+m_l)e_2-2\beta_l e_4} & {(1+m_{k+1})e_1 + \beta_l e_3}  & -3\beta_l e_2 & \beta_l e_1 \\
{(1+m_{k-1})e_1+3\beta_l e_3} & 0 & 0 & 0  \\
-3\beta_l e_2 & 0  & 0 & 0 \\
\beta_l e_1 & 0& 0 & 0 \\
\end{array}\right]^{(l)} 
 \left[ \begin{array}{rr}
0&e_1  \\
0&e_2  \\
0&e_3  \\
\beta_l^{-1}e_1&e_4  \\
\end{array}\right]^{(l+1)} \cdot\\
& \qquad \qquad \qquad \qquad \qquad \qquad \qquad \qquad \qquad \qquad \qquad \qquad \qquad \qquad \qquad \qquad   \left[\begin{array}{rr}e_1&m_{l+2}e_2\\0&e_1\end{array}\right]^{(l+2)}\cdots \left[ \begin{array}{r}m_d e_2\\e_1\end{array}\right]^{(d)}
\nonumber.
\end{align}
This explicit coefficient decomposition enables us to compute the ranks required to represent the FPUT equations via the different formats sketched in Fig.~\ref{fig:multiTT}.

\begin{rem}[Representation of the FPUT model with selection tensor $S$]\label{rem:FPUselection}
Each variable $x_k$ in the above decomposition of $\theta_l$ is represented by four different TT cores.
The FPUT model in the variant \eqref{eq:FPUMean} can thus be represented by a selection tensor $S$ with $n=4$ (see Fig.~\ref{fig:multiTT}b). The respective TT ranks are $r_k=4$.
\end{rem}

\begin{rem}[Representation of the FPUT model without selection tensor $S$]\label{rem:FPUwithoutselection}
Representing the governing equation of the FPUT model by a single TT (see Fig.~\ref{fig:multiTT}c) amounts to the decomposition of the tensor 
\begin{align*}
    \theta=\sum_{l=1}^d \theta_l\otimes e_l \, \in \mathbb{R}^{4^{\times d}\times d} \, .
\end{align*}
If $m_l=0$ and $\beta_l \neq 0$ the separation rank with respect to $\mathcal{P}_k=\{\{x_1,\ldots,x_k \}, \{x_{k+1},\ldots,x_{d},l\}\}$ is given by
\begin{align*}
r_k=
   \begin{cases} 4 &\text{if} \,\, k = 1\\
         4+k &\text{if} \,\, 1<k<d \\
         d &\text{if} \,\, k=d.
\end{cases}
\end{align*}
\end{rem}

Although the variants of the FPUT equation with $m_l=0$ are examples of a one-dimensional local interacting system with interaction range $(s_1,s_2)=(1,1)$ and separation rank $N=7$, we found representations with smaller ranks than estimated by Thm.~\ref{the:NN}. This is due to the special structure involving only nearest neighbor interactions (see Observation \ref{rem:lowrankpolynomials}).

\subsection{Randomized local interaction model}\label{app:extend}

In order to study generic models of one-dimensional local interacting systems, we focus here on randomly generated instances of systems with interaction range $(s_1,s_2)$.
Taking the first four Legendre polynomials as basis functions (see App.~\ref{app:L2orth}), we take each equation $f_l$ in a system $[f_1\ldots f_d]$ to be of the structure
\begin{equation}\label{eq:interactingModel}
\tilde{f}_l(x_{l-1},x_l,x_{l+1}) = \sum_{i_{l-1},i_l,i_{l+1}=1}^4 c_{i_{l-1} i_{l}i_{l+1}}\psi_{i_{l-1}}(x_{l-1})\psi_{i_l}({x_l})\psi_{i_{l+1}}(x_{l+1}) \,
\end{equation}
where $c_{i_{l-1} i_{l}i_{l+1}}$ are random coefficients, uniformly drawn on $[-1,1]$ and we set $x_0=x_{d+1}=0$. 
By construction, and using the notation of App.~\ref{app:fpu}, we note that each equation can be represented by a tensor $\tilde{\theta}_l$ with TT decomposition 
\begin{align*}
    \tilde{\theta}_l = A^{l-1} \cdot A^l \cdot A^{l+1}\, ,
\end{align*}
where we can choose core tensors with ranks $4$ as
\begin{align}
\nonumber
\small
A^{l-1} &=  \left[ \begin{array}{rrrr}
e_1 & e_2 & e_3 & e_4 
\end{array}\right]^{(l-1)} ,\\
A^l_{i_{l-1}i_{l+1}} &= \sum_{i_{l}=1}^4 c_{i_{l-1} i_{l}i_{l+1}} e^{(l)}_{i_l},\\
\nonumber
A^{l+1} &= \left[ \begin{array}{r}
e_1  \\
e_2  \\
e_3  \\
e_4  \\
\end{array}\right]^{(l+1)} .
\end{align}
We further note that the number of nonzero coefficients $c_{i_{l-1} i_{l}i_{l+1}}$ is an upper bound for the $\mathrm{CP}$-rank of $\tilde{\theta}_l$, and thus for the separation rank $N$ of the generated system (see Def.~\ref{def:nn}). 

\section{Extended numerical results}\label{app:numdiscuss}

\begin{figure}[tb]
    \centering
    \resizebox{0.7\linewidth}{!}{\input{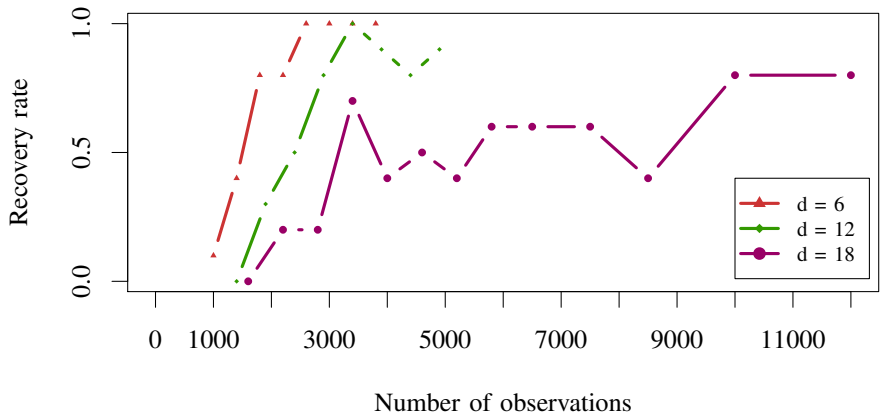}}
    \caption{The recovery rate (relative error threshold $10^{-6}$) for different number of observations $m$ and dimension $d$ for the FPUT equation with random $\beta_l$ and random mean field parameters $m_l$ (see \ref{app:fpu}). Each point is the average value of 10 trials.}
    \label{fig:exFPUTnl}
\end{figure}

\begin{figure}[tb]
    \centering
    \resizebox{0.7\linewidth}{!}{\input{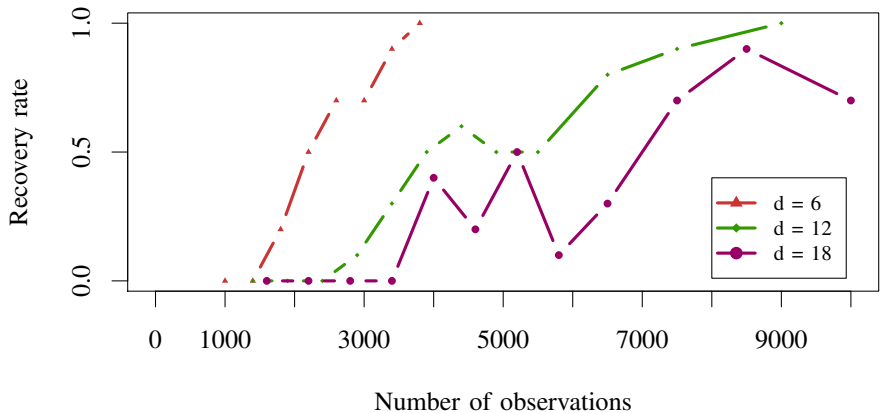}}
    \caption{The recovery rate (relative error threshold $10^{-6}$) for different numbers of observations $m$ and dimensions $d$ for random local interaction models trained with norm regularized ALS for 20 iterations (see \ref{app:extend}). Each point is the average of $10$ trials.}
    \label{fig:exFPUTrandomnl}
\end{figure}

In this section we present and discuss further numerical results.
While
Fig.~\ref{fig:stFPUT} shows the recovery rate for the FPUT equations 
\eqref{eq:FPUMean} 
with fixed constant $\beta_l = 0.7$ and $m_l = 0$, we will in this section turn our attention to more generic instances.
To this end, we draw the coefficients $\beta_l$ and $m_l$ i.i.d.\ uniformly at random from the interval $[-1,1]$. 
As explained in App.~\ref{app:fpu} the equations can still be represented by a TT model with rank $4$ using a selection tensor. 
We use this model trained by the norm-regularized ALS algorithm for the reconstruction (see also App.\ref{subsec:normreg}).
The resulting recovery rates averaged over $10$ trials for different number of observations $m$ are depicted in Fig.~\ref{fig:exFPUTnl}.

Going beyond the FPUT equation, the rank $4$ TT model with selection tensor also allows us to describe local interacting models in the form of Eq. \eqref{eq:interactingModel}, and as such we proceed to explore this more general setting.
To do so, we choose to use the first four Legendre polynomials (see App.~\ref{app:L2orth}) for the local function set $\{\psi_i\}$. 
We then draw instances of these equations that have $20$ non-vanishing coefficients $c_{ijk}$ per equation, with support selected independently uniformly at random for each equation. 
The values of the non-trivial coefficients are drawn i.i.d.\ uniformly at random from the interval $[-1,1]$. 
Using a sparse support allows us to keep the norm of $\tilde f$ moderately small. 
Fig.~\ref{fig:exFPUTrandomnl} shows the recovery rate for instances of such locally interacting equations for norm-regularized ALS.  Each point is the average of $10$ trial runs.

If a sufficient number of observations is provided, we observe a recovery rate close to $1$ for both random equation types with small $d=6,12$ using ALS. 
The numerical results in the Fig.~\ref{fig:stFPUT}, \ref{fig:exFPUTnl} and \ref{fig:exFPUTrandomnl} additionally demonstrate that for larger $d$, here $d=18$, the recovery rates were poorly improving when increasing the number of observations $m$ in the norm regularized ALS. In cases, where recovery was not achieved, especially for larger number of observations $m$, the convergence behavior of the relative error and the residual of the iterate $\theta_k$ indicated that ALS got stuck in local minima. This is a known issue with alternating optimization schemes like ALS. In the two randomized models, see App.~\ref{app:fpu}~and~\ref{app:extend}, we have three random parts. Firstly, the governing equations have random coefficients. Secondly, the observations are drawn randomly, and thirdly the ALS like recovery schemes are initialized randomly. 
In order to test the intuition that successful recovery depends strongly on the random initialization we ran an additional simple experiment: 
If the error after $25$ iterations of the regularized ALS was not less than $10^{-6}$, we restarted the method, up to $5$ times, with a different random initialization. 
For different numbers $m$ of observations and $d=18$ we get for the restarted version the recovery rates depicted in Table \ref{tab:restartreg}. We also state the averaged number of restarts, which is $4$ in case of no success of the recovery method.

\begin{table}[tb]
\centering
\begin{tabular}{|l|l|l|l|}
\hline
Dimension $d$ & Number of observations $m$ & Recovery Rate & Averaged Number of Restarts \\
\hline\hline
18 & 1000 & 0 of 10 & 4\\\hline
18 & 2000 & 1 of 10 & 3.8\\\hline
18 & 3000 & 6 of 10  & 2.5 \\\hline
18 & 4000 & 10 of 10 & 1.3\\\hline
18 & 5000 &  10 of 10 & 0.2\\\hline
18 & 6000 &  10 of 10 & 0.3\\\hline
18 & 7000 &  10 of 10 & 0.4\\\hline
\end{tabular}
\vspace{.5cm}
\caption{Recovery rates for the restarted regularized ALS for randomized governing equations as in App.~\ref{app:extend}. Again, recovery is achieved, if the relative error drops below $10^{-6}$. The last column states the number of restarts (max. $5$).
\label{tab:restartreg}}
\end{table}

\section{$L_2([-1,1])$-orthogonal polynomials} \label{app:L2orth}
As univariate basis functions we have used in our simulations the Legendre polynomials
\begin{align*}
\psi_1(x) &=  1,\\
\psi_2(x) &= x,\\
\psi_3(x) &= \frac{1}{2}(3x^2-1),\\
\psi_4(x) &=\frac{1}{2}(5x^3-3x).\\
\end{align*}
The coefficients are chosen in such a way that the polynomials are $L_2$ orthogonal on $[-1,1]$.

\section{Details on the numerical implementation} \label{app:algorithm}

In this section we provide a detailed discussion of the numerical methods, ALS and SALSA, which were applied in our numerical experiments.
We furthermore derive the complexity of the algorithms, which enables us to compare the computational requirements necessary for applying these algorithms with the different tensor models we have introduced (see Section \ref{sec:examples}).
All code and numerical examples are openly available at the associated GitHub repository \cite{anon-git-repo} (\url{https://github.com/RoteKekse/systemrecovery}).

\subsection{Implementational notes on ALS} \label{subsec:normreg}
In the examples we have chosen to study, the selection tensor has a non-trivial low-rank kernel, and as such for the tensor network model involving the selection tensor problems can arise in the application of ALS due the potential for an increasing norm.
In order to address this issue, 
we regularize the problem $(P_k)$ by introducing a norm penalty in the following way
\begin{equation} 
\min_{{A}^k} \| f (\tilde{A}^1,\dots,{A}^k,\dots,\tilde{A}^d) - y \|_F^2 + \lambda \|{A}^k\|_F^2 \,.\tag{$P_k^{\mathrm{reg}}$}
\end{equation}
A core $A^k$ optimizing the loss in $(P_k^{\mathrm{reg}})$ is found by the first order condition (see Fig.~\ref{fig:ALS}):
\begin{equation}\label{eq:regals}
    \begin{split}
        0&=\nabla_{A^k}\| f (\tilde{A}^1,\dots,A^k,\dots,\tilde{A}^d) - y \|_F^2+\lambda \nabla_{A^k}  \|{A}^k\|_F^2\\
        &= 2 \big[ \braket{A^k\nabla_{A^k} f, \nabla_{A^k} f } - \braket{y, \nabla_{A^k} f } + \lambda A_k
        \big],
    \end{split}
\end{equation}
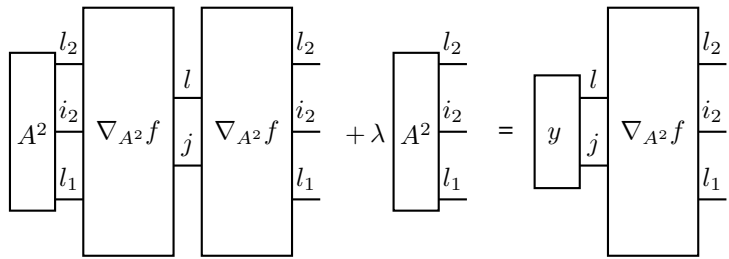
\begin{figure}[tb]
    \centering
   \begin{tikzpicture}[scale=0.3,thick]
   \draw (-7.5,1) rectangle (-3.5,-10); 
   \coordinate[label=above:{$ \nabla_{A^2} f$}] (A) at (-5.5,-5.5);
   \draw (-3.5,-6)--(-2.25,-6) node[midway,above] {$j$};
   \draw (-3.5,-3)--(-2.25,-3) node[midway,above] {$l$};
   \draw (-7.5,-4.5)--(-8.75,-4.5) node[midway,above] {$i_2$};
   \draw (-7.5,-1.5)--(-8.75,-1.5) node[midway,above] {$l_2$};
   \draw (-7.5,-7.5)--(-8.75,-7.5) node[midway,above] {$l_{1}$};
   \draw (-8.75,-8) rectangle (-10.75,-1);
   \coordinate[label=above:{$ {A}^2 $}] (A) at (-9.75,-5.5);

\begin{scope}[yscale=1,xscale=-1,shift={(5.75,0)}]
   \draw (-7.5,1) rectangle (-3.5,-10); 
   \coordinate[label=above:{$ \nabla_{A^2} f$}] (A) at (-5.5,-5.5);
   \draw (-7.5,-4.5)--(-8.75,-4.5) node[midway,above] {$i_2$};
   \draw (-7.5,-1.5)--(-8.75,-1.5) node[midway,above] {$l_2$};
   \draw (-7.5,-7.5)--(-8.75,-7.5) node[midway,above] {$l_{1}$};
\end{scope}

\begin{scope}[yscale=1,xscale=-1,shift={(-12.25,0)}]
   \draw (-7.5,1) rectangle (-3.5,-10); 
   \coordinate[label=above:{$ \nabla_{A^2} f$}] (A) at (-5.5,-5.5);
   \draw (-3.5,-6)--(-2.25,-6) node[midway,above] {$j$};
   \draw (-3.5,-3)--(-2.25,-3) node[midway,above] {$l$};
   \draw (-2.25,-7) rectangle (-0.25,-2);
   \coordinate[label=above:{$y$}] (A) at (-1.125,-5.5);
   \draw (-7.5,-4.5)--(-8.75,-4.5) node[midway,above] {$i_2$};
   \draw (-7.5,-1.5)--(-8.75,-1.5) node[midway,above] {$l_2$};
   \draw (-7.5,-7.5)--(-8.75,-7.5) node[midway,above] {$l_{1}$};
\end{scope}

\begin{scope}[yscale=1,xscale=1,shift={(17,0)}]
   \coordinate[label=above:{$ + \, \lambda $}] (A) at (-12,-5.5);
   \draw (-8.75,-8) rectangle (-10.75,-1);
   \coordinate[label=above:{$ {A}^2 $}] (A) at (-9.75,-5.5);
   \draw (-7.5,-4.5)--(-8.75,-4.5) node[midway,above] {$i_2$};
   \draw (-7.5,-1.5)--(-8.75,-1.5) node[midway,above] {$l_2$};
   \draw (-7.5,-7.5)--(-8.75,-7.5) node[midway,above] {$l_{1}$};
\end{scope}

\coordinate[label=left:{=}] (T) at (12,-4.5);

\end{tikzpicture}
    \caption{Least squares optimization of the tensor train component $A^2$, by the solution of a linear equation,
    the first order minimum condition of $(P_k)$.}
    \label{fig:ALS}
\end{figure}
which is just the solution of a linear equation system.
We thereby used the scalar product corresponding to the Frobenius norm.
In our numerical experiments, the results of which are seen in Figures ~\ref{fig:stFPUT}, \ref{fig:exFPUTnl} and \ref{fig:exFPUTrandomnl}, we initialized $\lambda$ at $\lambda = 1$, and after each optimization sweep through all cores modified $\lambda$ via the simple heuristic ${\lambda_{\rm  new} = \lambda_{\rm old} / 10 }$.
In the restarted version of the algorithm ( described in App.~\ref{app:numdiscuss} and Table~\ref{tab:restartreg}) we 
have used the following heuristic
\begin{align}\label{eq:regupdate}
\lambda_{\rm  new} &= \min \left\{0.1 \frac{\| f (\tilde{A}^1,\dots,{A}^k,\dots,\tilde{A}^d) - y \|_F^2}{\|y\|_F\|{A}^k\|_F}, \frac{\lambda_{\rm old}}{4}\right\},
\end{align}
which aims at balancing the two terms in $(P_k^{reg})$.\\
For an efficient computation of the gradient $\nabla_{A^k}f$ (see Fig.~\ref{fig:gradient}), we use so called stacks (see \citet{wolf_low_rank_2019}). The benefit are achieved by a trade-off between memory and computational demand. We store the dictionary tensor as a list of $d$ $m\times \tilde p$ matrices, therefore having a storage cost of $dm \tilde p$.
For the examples in this work, the selection tensor can also be stored by a list of $d$ $n\times d$ matrices with storage consumption $d^2n$. 
The stacks used in the above numerical experiments, for one update step, then look as follows (here for the left stack):
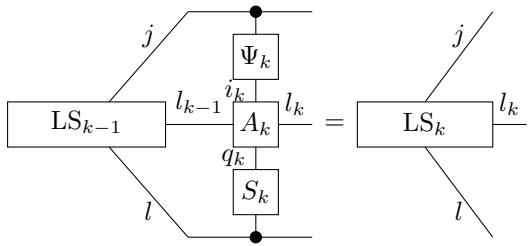
\begin{center}
\begin{tikzpicture}[scale=0.3]
    \draw (-10,-1) rectangle (-3,1);
     \draw (-5.5,1)--(-2,5) node[midway,above] {$j$};
     \draw (-5.5,-1)--(-2,-5) node[midway,below] {$l$}; 
     \draw (-3,0)--(0,0) node[midway,above] {$l_{k-1}$};
    \coordinate[label=above:{$\mathrm{LS}_{k-1}$}] (A) at (-6.5,-0.9);

    \draw (0,-1) rectangle (2,1);
    \draw (0,2) rectangle (2,4);
    \draw (0,-4) rectangle (2,-2);
    \coordinate[label=above:{$A_k$}] (A) at (1,-0.9);
    \coordinate[label=above:{$\Psi_k$}] (A) at (1,2.1);
    \coordinate[label=above:{$S_k$}] (A) at (1,-3.9);

    \draw (2,0)--(3.5,0) node[midway,above] {$l_k$};
    \draw (1,1)--(1,2) node[midway,left] {$i_k$};
    \draw (1,-1)--(1,-2) node[midway,left] {$q_k$};
    \draw (-2,5)--(3.5,5) ;
    \draw (-2,-5)--(3.5,-5) ;
    \draw (1,4)--(1,5) ;
    \draw (1,-4)--(1,-5) ;
    \draw[fill=black] (1,-5) circle (0.25cm);
    \draw[fill=black] (1,5) circle (0.25cm);

    \node at (4.5,0) {$=$};

    \begin{scope}[shift={(0.5,0)}]
\draw (5,-1) rectangle (11,1);
     \draw (8,1)--(11,5) node[midway,above] {$j$};
     \draw (8,-1)--(11,-5) node[midway,below] {$l$}; 
     \draw (12.5,0)--(11,0) node[midway,above] {$l_k$};
    \coordinate[label=above:{$\mathrm{LS}_{k}$}] (A) at (8,-0.9);
    \end{scope}
\end{tikzpicture}
\end{center} 
One such stack update has complexity of $\mathcal{O}(md\tilde{p}r^2n)$ for the Format (\ref{fig:multiTT}b), where $m$ is the number of observations, $\tilde{p}$ is the number of univariate basis functions, $d$ is the number of variables and number of equations,
$r$ is the TT-rank  
of the given format and $n$ is the number of activation patterns.
 If one can store the left and right stacks building the local linear 
 equation systems $(P_k)$ (see also Fig. \ref{fig:ALS})  ALS schemes are much more efficient. Furthermore, the aforementioned local linear equation system's solution $A_k$ of \eqref{eq:regals} has dimensions $r_{k-1}r_kn_k\tilde{p}$ for the selection tensor format (\ref{fig:multiTT}b). Building the local linear operator has complexity  $\mathcal{O}(md(\tilde{p}r^2n)^2)$. The complexity to solve the system of linear equations, using standard linear equation solvers, e.g. Gauss algorithm, is $\mathcal{O}((r^2n\tilde{p})^3)$. We therefore conclude that the contraction building the local linear operator is the most expensive part, if we can bound $n$ and $r$.

\subsection{Implementational notes on SALSA} \label{subsec:appsalsa}
The SALSA scheme, as introduced by \citet{grasedyck_stable_2019}, iteratively performs a rank-adaptive variant of the update step $(P_{k-s})$ (see Section~\ref{sec:optimization}).
The parameters $\omega$ and $\epsilon$ are adjusted during the iterations, where we state in the following our choices for the initalization and update rule.
For the rank adaptive numerical experiment we have used the initialization ${\omega_{\rm start} = 1}$ and $\epsilon_{\rm start}=0.2$, furthermore we set ${r_{\rm min} = 2}$, ${s_{\rm min} = 0.2}$, ${\omega_{\rm min} = 1.05}$,  ${r_{\rm start} = (1,\ldots,1)}$ and ${c = 0.01}$. The parameters are updated in the following algorithm:

\begin{algorithm}[H]{Stabilized alternating least squares approximation (SALSA)}\\
Input: selection tensor $S$, dictionary tensor $\Phi$, right hand side $y$\\ 
Output: iterate solution $\theta$
\begin{itemize}
    \item[1)] Fix $r_{\rm min} \in\mathbb{N}$, $r = r_{\rm  start}\in\mathbb{N}^{d-1}$, $\epsilon= \epsilon_{\rm start},\omega = \omega_{\rm start} > 0$.
    \item[2)] Initialize randomly ${\theta= A^1\cdot \ldots \cdot A^d}$, ${A^k \in\mathbb{R}^{(1+r_{\rm min})\times \tilde{p} \times (1+r_{\rm min})}}$
    \item[3)] Solve ($P_{k-\mathrm{s}}$) for all $k = 1 ,\dots, d$ at least once for ${A_k \in\mathbb{R}^{(r_{k-1}+r_{\rm min})\times \tilde{p}\times (r_k+r_{\rm min})}}$.
    \item[4)] For $k = 1,\ldots ,d-1$ set the new rank $r_k^{(new)}$ to the number of singular values $\sigma_{k,j}, k = 1,\ldots ,d-1, j = 1,\dots, r_k+r_{\rm min}$ which are greater than $\epsilon$.
    \item [5)] If the rank $r_k$ increased in step 4) add new singular values of size $c\epsilon, 0 < c < 1$ to the $k$th virtual index. This changes the tensor only little.
    \item [6)] Decrease $\omega$ and $\epsilon$. In our numerical experiments, we have used for ${R(\theta) := \|\Phi \theta - y \|_F}$
    \begin{align*}
       &\omega_{\rm new} = \min\{\sqrt{R(\theta)}, \frac{\omega_{\rm  old}}{\omega_{\rm min}}\},\\
       &\epsilon_{\rm new} = s_{\rm min}R(\theta), 
    \end{align*} 
    which heuristically showed the best performance.
    \item [7)] Repeat steps 3) to 6) until the residual $ R(\theta)$ is smaller than some prescribed threshold.
\end{itemize}
\end{algorithm}
Again, for an efficient implementation of the network contraction we can employ stacks. For the single TT format, similarly to the ALS case, the complexity of the stack updates is in $\mathcal{O}(m\tilde{p}dr^2)$ and building the local linear operator for the system of linear equations (see \ref{fig:ALS}) is $\mathcal{O}(md(\tilde{p}r^2)^2)$. The computation of the solution of the linear equation system has a demand in $\mathcal{O}((r^2\tilde{p})^3)$. 

We have used SALSA for the single TT format~(\ref{fig:multiTT}c) and by Theorem~\ref{the:NN} the ranks increase linearly, while in the selection Format~(\ref{fig:multiTT}b) the ranks are bounded. 
It follows that for the example systems we study, the scaling of the update steps for the single TT formal is unfavourable for large $d$, which can be avoided by use of the selection tensor format.

\end{document}